\documentclass[a4paper]{amsart}
\usepackage{amssymb}
\usepackage{amsmath}
\usepackage{xypic}

\newtheorem{thm}{Theorem}[section]
\newtheorem*{thm*}{Theorem}
\newtheorem{thmintro}{Theorem}

\newtheorem{cor}[thm]{Corollary}
\newtheorem{defn}[thm]{Definition}
\newtheorem{ex}[thm]{Example}
\newtheorem{lemma}[thm]{Lemma}

\newtheorem{prop}[thm]{Proposition}

\newtheorem{rem}[thm]{Remark}

\newcommand{\ZZ}{{\mathbb  Z}}
\newcommand{\RR}{{\mathbb  R}}
\newcommand{\CC}{{\mathbb  C}}

\renewcommand{\cH}{{\mathcal H}}

\renewcommand{\cH}{{\mathcal H}}

\newcommand{\E}{{\mathcal E}}
\newcommand{\oE}{\overline{{\mathcal E}}}

\newcommand{\F}{{\mathcal F}}
\newcommand{\oF}{\overline{{\mathcal F}}}

\newcommand{\mR}{{\mathbb R}}

\newcommand{\oL}{\overline{L}}
\newcommand{\Oq}{{\mathbf O(q)}}

\newcommand{\soq}{{{\mathfrak{so}}(q)} }
\newcommand{\SOq}{{\mathbf{SO}(q)}}
\newcommand{\bas}{{\rm bas}\, }

\newcommand{\whL}{{\widehat L}}
\newcommand{\whM}{{\widehat M}}

\newcommand{\mfg}{\mathfrak{g}}
\newcommand{\mfh}{\mathfrak{h}}
\newcommand{\mfb}{\mathfrak{b}}
\newcommand{\mfa}{\mathfrak{a}}
\newcommand{\mfs}{\mathfrak{s}}
\newcommand{\mfk}{\mathfrak{k}}
\newcommand{\mft}{\mathfrak{t}}
\newcommand{\mfm}{\mathfrak{m}}

\newcommand{\umfg}{\underline{\mathfrak{g}}}

\newcommand{\Ad}{\operatorname{Ad}}

\newcommand{\Ann}{\operatorname{Ann}}

\newcommand{\codim}{\operatorname{codim}}
\newcommand{\depth}{\operatorname{depth}}
\newcommand{\grade}{\operatorname{grade}}
\newcommand{\rk}{\operatorname{rank}}
\newcommand{\trdim}{\operatorname{trdim}}
\newcommand{\Hess}{\operatorname{Hess}}

\numberwithin{equation}{section}

\begin{document}

\title[Equivariant Basic Cohomology]{Equivariant Basic Cohomology of Riemannian Foliations}

\author{Oliver Goertsches}
\address{Oliver Goertsches, Mathematisches Institut, Universit\"at zu K\"oln, Weyertal 86-90, 50931 K\"oln, Germany}
\email{ogoertsc@math.uni-koeln.de}
\author{Dirk T\"oben}
 \thanks{DT was supported  by the Schwerpunktprogramm SPP 1154  of the DFG}
\address{Dirk T\"oben, Mathematisches Institut, Universit\"at zu K\"oln, Weyertal 86-90, 50931 K\"oln, Germany}
\email{dtoeben@math.uni-koeln.de}

\date{}


\begin{abstract}
The basic cohomology of a Riemannian foliation on a complete manifold with all leaves closed is the cohomology of the leaf space. In this paper we introduce various methods to compute the basic cohomology in the presence of both closed and non-closed leaves in the simply-connected case (or more generally for Killing foliations): We show that the total basic Betti number of the union $C$ of the closed leaves is smaller than or equal to the total basic Betti number of the foliated manifold, and we give sufficient conditions for equality. If there is a basic Morse-Bott function with critical set equal to $C$ we can compute the basic cohomology explicitly. Another case in which the basic cohomology can be determined is if the space of leaf closures is a simple, convex polytope.
Our results are based on Molino's observation that the existence of non-closed leaves yields a distinguished transverse action on the foliated manifold with fixed point set $C$. We introduce equivariant basic cohomology of transverse actions in analogy to equivariant cohomology of Lie group actions enabling us to transfer many results from the theory of Lie group actions to Riemannian foliations. The prominent role of the fixed point set in the theory of torus actions explains the relevance of the set $C$ in the basic setting.
\end{abstract}

 \maketitle

\tableofcontents

\section{Introduction} \label{sec-intro}
One of the fundamental topological invariants of a foliation are its basic Betti numbers introduced by Reinhart \cite{Reinhart}. In this paper we want to show how dynamical aspects of a Riemannian foliation are linked to these numbers.
A Riemannian foliation $\F$ is a foliation by locally equidistant leaves (see Section \ref{subsec:strkillalg}). Examples are foliations given by the connected fibers of Riemannian submersions or the orbit decomposition of locally free isometric actions. In the first example the leaves have trivial holonomy and are closed, which is not the case for general Riemannian foliations. The partition $\oF$ by leaf closures of $\F$, a so-called singular Riemannian foliation, therefore usually differs from $\F$; this difference reflects the dynamical behavior of $\F$.
Let us look at an example of the second kind illustrating this behavior. For rationally independent reals $s_1,s_2$ we consider the action $\RR\times S^3\to S^3;(t,(z_1,z_2))\mapsto (e^{2\pi is_1t}z_1,e^{2\pi is_2t}z_2)$. The orbit decomposition defines a one-dimensional Riemannian foliation which we call the Hopf flow for further reference. It has exactly two closed leaves. One of the motivating questions for us was whether there is a topological relation between the union of the closed leaves and the foliated space, just as by the Poincar\'e-Hopf Theorem and a Theorem by Bott \cite{Bott} there is one between the zeroes of a vector field/Killing field and the (Riemannian) manifold.  In fact we will see that this is indeed the case; the role of the vector field/Killing field will be taken by an algebra of transverse Killing fields which occurs naturally in our situation.

By Molino's structure theory for Riemannian foliations on complete manifolds \cite{Molino} (see Theorems \ref{thm:Killing} and \ref{thm:equivariantMolino} in this paper), the leaf closures are the orbits of leaves under the action of an abelian Lie algebra $\mfa$ of transverse Killing fields, called the structural Killing algebra of $(M,\F)$, at least if $M$ is simply-connected or more generally for Killing foliations (Definition \ref{defn:Killing}).  We will speak of an isometric transverse action of $\mfa$ on $(M,\F)$.

The idea now is to pass from the basic cohomology $H^*(M,\F)$ whose Betti numbers we are interested in to $\mfa$-equivariant basic cohomology $H^*_\mfa(M,\F)$ which is defined in analogy to ordinary equivariant cohomology via the language of $\mfg^\star$-algebras. $H^*_\mfa(M,\F)$ has a natural algebra structure over the polynomial ring $S(\mfa^*)$. Just as for torus actions, a special class of transverse actions will turn out to be of particular importance, namely those that are equivariantly formal: these are actions for which $H^*_\mfa(M,\F)\cong S(\mfa^*)\otimes H^*(M,\F)$ as $S(\mfa^*)$-modules. Many results from the theory of torus actions obtained via equivariant cohomology can be transferred to Riemannian foliations. 
For instance, we have the following (Theorem \ref{thm:boreleq}) in spirit of the above motivation  -- in the subsequent theorems, assume that $\F$ is a transversely oriented Killing foliation of codimension $q$ on a complete manifold $M$ with structural Killing algebra $\mfa$ such that $M/\oF$ is compact, and denote by $C$ the union of the closed leaves of $\F$.
\begin{thmintro}\label{intro:Borel} We have 
$$
\dim H^*(C/\F)=\dim H^*(C,\F)\leq \dim H^*(M,\F),
$$
and equality is valid if and only if the $\mfa$-action on $(M,\F)$ is equivariantly formal.
\end{thmintro}
In case of an equivariantly formal action the Poincar\'e series $P_t(M,\F)$ of $H^*(M,\F)$ can be recovered from the equivariant Poincar\'e series $P_t^\mfa(M,\F)$ of $H^*_\mfa(M,\F)$. We find verifiable sufficient conditions for equivariant formality which at the same time guarantee the computability of $P_t^\mfa(M,\F)$; hence, these conditions enable us to compute $P_t(M,\F)$. In this way we prove (Theorem \ref{thm:basicperfect})
\begin{thmintro}\label{intro:MorseBott}
 Any basic Morse-Bott function  $f:M\to\RR$ with critical set $C$ is perfect, i.e.,
$$
P_t(M,\F)=\sum_N t^{\lambda_N} P_t(N/\F),
$$
where $N$ runs over the components of $C$, and $\lambda_N$ is the index of $f$ at the critical manifold $N$.
\end{thmintro}
and (Theorem \ref{thm:polytopebasicps})
\begin{thmintro}\label{intro:polytope} 
If $M/\oF$ is a simple, convex polytope of dimension $\dim \mfa$ such that the vertices correspond to the closed leaves, then 
$$
P_t(M,\F)=\sum_i \lambda_i t^{q-2i}(1-t^2)^i,
$$
where $\lambda_i$ denotes the number of faces of $M/\oF$ of dimension $i$.
\end{thmintro}

The latter theorem follows from (Corollary \ref{cor:basicps})

\begin{thmintro}\label{intro:series}
If the $\mfa$-action on $(M,\F)$ is equivariantly formal, then
$$
P_t(M,\F)=\sum_{X} {t^{\codim X }}{(1-t^2)^{{\rm trdim\,}L_X}}P_t(X/\oF),
$$
where $X$ runs over the components of infinitesimal orbit type manifolds and $L_X$ is any leaf in $X$. Here, the transverse dimension of the leaf $L_X$ is defined as ${\rm trdim\,}L_X=\dim \oL_X-\dim L_X$.
\end{thmintro}

To prove these results, two techniques are equally important: On the one hand, we transfer concepts and results from ordinary equivariant cohomology to the equivariant basic case, like a Borel localization theorem (Theorem \ref{thm:borelloc}):

\begin{thmintro}\label{intro:Localization} The natural restriction map
$$
\widehat H^*_\mfa(M,\F)\to \widehat H^*_\mfa(C,\F)
$$
is an isomorphism, where $\widehat H$ denotes localization at the zero ideal in $S(\mfa^*)$.
\end{thmintro}
On the other hand, we establish a correspondence principle between a transversely oriented Riemannian foliation $\F$ and an associated $\SOq$-manifold $W$ using Molino theory. It allows us to translate basic data of the foliation to equivariant data of the $\SOq$-action and back. Most importantly it yields (Proposition \ref{prop:commutingaction})
\begin{thmintro}\label{intro:correspondence}
$H^*_\mfa(M,\F)\cong H^*_\SOq(W)$ as graded rings.
\end{thmintro}
This links equivariant basic cohomology of $(M,\F)$ closely to ordinary equivariant cohomology of the associated $\SOq$-manifold $W$. In particular (Proposition \ref{prop:eqformalCM}), it relates the condition of equivariant formality of the transverse $\mfa$-action to the Cohen-Macaulay property of the $\SOq$-action on $W$. As all isotropy algebras of the $\SOq$-action are abelian, the latter property behaves like in the case of torus actions, for which it was investigated in \cite{GT2009}.

Besides being useful as a tool to prove results on classical basic cohomology like the theorems listed above, we think that equivariant basic cohomology is interesting in its own right as a new invariant for Killing foliations.  Results in this direction are the  calculations in Section \ref{sec:basiccalculations}, the computability of the equivariant basic Poincar\'{e} series in terms of the stratification of $M$ induced by the $\mfa$-action even in case the $\mfa$-action is not equivariantly formal (Theorem \ref{thm:basicpoincareseries}) or the equivariant basic Morse-Bott inequalities (Theorem \ref{thm:basicMorseBott}). It is conceivable that our technique allows one to translate more results on torus actions and their equivariant cohomology to the theory of Killing foliations.   

In Section \ref{sec:transverseaction} we define transverse actions on foliated manifolds and isometric transverse actions for Riemannian foliations. We mention Molino's result to be explained in detail in Section \ref{sec:correspondence} that a Riemannian foliation on a simply-connected manifold, or more generally a Killing foliation, has a natural isometric transverse action of the so-called structural Killing algebra whose orbits are the leaf closures.

We observe in Section \ref{sec:equivbasiccoh} that a transverse action of a Lie algebra $\mfg$ on a foliated manifold $(M,\F)$ induces the structure of a $\mfg^\star$-algebra on  the complex of basic forms  $\Omega(M,\F)$. This leads to the notion of equivariant basic cohomology denoted by $H^*_\mfg(M,\F)$. In the remainder of the section we prove elementary properties of this equivariant cohomology theory and introduce the notion of equivariant formality.

In Section \ref{sec:correspondence} we establish the correspondence between $(M,\F)$ and the $\SOq$-manifold $W$ described above, and in particular Theorem \ref{intro:correspondence}. Besides equivariant cohomology this correspondence also links the respective isotropy algebras and stratifications.

The Borel-type localization theorem (Theorem \ref{intro:Localization}) and its immediate corollaries (e.g., Theorem  \ref{intro:Borel}) are covered in Section \ref{sec:Borel}.

We use the correspondence principle to reduce (equivariant) basic Morse Bott theory on $(M,\F)$ to (equivariant) Morse-Bott theory on $(W,\SOq)$ in Section \ref{sec:basicMorseBott}. In particular we prove Theorem \ref{intro:MorseBott} and the equivariant basic Morse-Bott inequalities.

With the help of the equivariant Gysin sequence, which is investigated in Appendix \ref{sec:appA}, we calculate in Section \ref{sec:Poincareseries} for a certain class of $G$-manifolds $W$ the equivariant Poincar\'{e} series of $H^*_G(W)$. Via the correspondence principle, this yields an expression for the Poincar\'{e} series of $H^*_\mfa(M,\F)$ (Theorem \ref{thm:basicpoincareseries}), which in the special case of an equivariantly formal $\mfa$-action gives Theorem \ref{intro:series}.

The expression for the Poincar\'{e} series of $H^*_\mfa(M,\F)$ mentioned above contains the Poincar\'{e} polynomials of the spaces $X/\oF$, where $X$ runs through the components of infinitesimal orbit type manifolds. If $M/\oF$ is a simple, convex polytope as a stratified space, these polynomials are as simple as possible.  In Section \ref{sec:polytope} it turns out that in this case the $\mfa$-action is automatically equivariantly formal and we  derive Theorem \ref{intro:polytope}.

{\em Acknowledgements:} The second author is grateful to Steven Hurder for various discussions on the topology of foliations.

\section{Transverse actions}\label{sec:transverseaction}
Let $\F$ be a foliation on a manifold $M$. By $\Xi(\F)$ we denote the space of differentiable vector fields on $M$ that are tangent to the leaves. A vector field $X$ on $M$ is said to be {\em foliate} if for every $Y\in \Xi(\F)$ the Lie bracket $[X,Y]$ also belongs to $\Xi(\F)$. See \cite[Proposition 2.2]{Molino} for a geometric characterization of foliate fields. The set $L(M,\F)$ of foliate fields is the normalizer of $\Xi(\F)$ in the Lie algebra $\Xi(M)$ of vector fields on $M$ and therefore a Lie sub-algebra of $\Xi(M)$. We call the projection of a foliate field $X$ to $TM/T\F$ a {\em transverse} field. The set $l(M,\F)=L(M,\F)/\Xi(\F)$ of transverse fields is also a Lie algebra inheriting the Lie bracket from $L(M,\F)$.
\begin{defn}\label{defn:transverseaction}
A  {\em transverse action} of a finite-dimensional Lie algebra $\mfg$ on the foliated manifold $(M,\F)$ is a Lie algebra homomorphism $\mfg\to l(M,\F)$. 
\end{defn}
Given a transverse action of $\mfg$, we will denote the transverse field associated to $X\in \mfg$ by $X^*\in l(M,\F)$.
If $\F$ is the trivial foliation by points, this notion coincides with the usual notion of an infinitesimal action on the manifold $M$.

We assume our actions to be effective unless otherwise stated, i.e., $X=0$ whenever $X^*=0$. Let $U$ be a foliation chart domain and $p:U\to U/\F$ the projection. Then the transverse action of $\mfg$ descends to an ordinary action of $\mfg$ on $U/\F$. In this way $p$ becomes $\mfg$-equivariant. 
We consider the distribution $\mfg \cdot \F$ of varying rank defined by $(\mfg\cdot \F)_p=\{X^*_p\mid X\in \mfg\} + T_p\F$. Note that $\{X^*_p\mid X\in \mfg\}$ by definition is not a subspace of $T_pM$; nevertheless, $(\mfg\cdot \F)_p$ is well-defined. The distribution $\mfg\cdot\F$ is spanned by the set of vector fields $D=\Xi(\F)\cup \{Y\in \Xi(M)\mid \mbox{for some}\ X\in\mfg\ \mbox{the associated transverse field of}\ Y\ \mbox{coincides with}\ X^*\}$. One can show that $\mfg\cdot\F$ is invariant under the local flows of vector fields in $D$ (it follows $\mfg\cdot\F=P_D$ in the notation of Sussmann \cite{Sussmann}). By \cite[Theorem 4.2]{Sussmann}, $\mfg\cdot\F$ has the maximal integral manifolds property \cite[p.~178]{Sussmann}. The orbits of $D$ as defined by Sussmann constitute these maximal integral manifolds \cite[Theorem 4.1]{Sussmann}, and coincide with the accessible sets of Stefan \cite{Stefan}. By  \cite[Theorem 4]{Stefan} these define a singular foliation which is everywhere tangent to $\mfg\cdot\F$. We will use $\mfg\cdot\F$ synonymously for the distribution and its singular foliation. The leaf of $\mfg\cdot \F$ through a point $p$ contains the leaf $L_p\in \F$ through $p$. We will denote it by $\mfg\cdot L_p$ and refer to it as the {\em $\mfg$-orbit of the leaf $L_p$}.

If the rank of $\mfg\cdot \F$ is constant, then $\mfg\cdot\F$ is a foliation.  If for any $v\in TM/T\F$ there is an $X\in\mfg$ with $X_p^*=v$, where $p$ is the foot point of $v$, we call the action {\em transitive}; then the orbit of a leaf is a connected component of $M$. The action is {\em free}, if $X=0$ whenever $X_p^*=0$ for some $p\in M$.

\subsection{Isometric transverse actions}\label{subsec:strkillalg} In this subsection we will recall the definition of Riemannian foliations and introduce the notion of an isometric transverse action. As the main example to be kept in mind for Section \ref{sec:equivbasiccoh} we mention a result of Molino, which in our language states the existence of a distinguished (local) transverse action on a manifold with a Riemannian foliation whose orbits are the leaf closures. The construction of this action will be explained in Section \ref{sec:Molinobundle}.

Let $\F$ be a foliation and $g$ a metric on $TM/T\F$, called {\em transverse metric}. 
Let $U\subset M$ be a foliation chart domain with projection $p:U\to V\subset \RR^q$ defining $\F$ on $U$ as the fibers. If for every such projection $p$ any two vectors $v,w\in TU/T\F$ with $p_*v=p_*w$ have the same length with respect to $g$ we say that $\F$ is a {\em Riemannian foliation}. In this case $V\cong U/\F$ can be endowed with a metric $p_*g$ such that $p:(U,g)\to (U/\F,p_*g)$ becomes a Riemannian submersion. In this sense a Riemannian foliation is locally given by Riemannian submersions. We pull back the Levi-Civita connection from $(U/\F,p_*g)$ to $U$. These connections coincide on overlaps because of the uniqueness of the Levi-Civita connection. We obtain a basic Riemannian connection, where basic means that the local transverse fields $l(U,\F)$ are parallel along the leaves with respect to this connection.

Let $\F$ be a Riemannian foliation with transverse metric $g$ on a complete manifold $M$. A vector field $X$ with $L_Xg=0$ is automatically foliate \cite[Lemma 3.5]{Molino} and we call its transverse field a {\em transverse Killing field}. The defining condition is closed under the Lie bracket and therefore the set of transverse Killing fields denoted by ${\mathfrak{iso}}(M,\F,g)\subset l(M,\F)$ is a Lie algebra. A transverse action $\mfg\to {\mathfrak{iso}}(M,F,g)$ is called an {\em isometric transverse action}. As above let $U$ be a foliate chart domain and $p:(U,g)\to (U/\F,p_*g)$ the projection. Then the transverse action of $\mfg$ descends to an action of $\mfg$ on $U/\F$ by Killing fields with respect to the projected metric $p_*g$. Since a Killing field $X$ on $U/\F$ is uniquely determined by $X_p$ and $(\nabla X)_p$ for any $p\in U/\F$, so is every transverse Killing field; here $\nabla$ denotes the Levi-Civita connection of $(U/\F,p_*g)$. Note that $\mfg\cdot\F$ is a singular Riemannian foliation (see \cite[p.~189]{Molino}) in the sense of Molino. The proof given in \cite[Proposition 6.2]{Molino} is stated for the partition $\oF=\{\oL\mid L\in \F\}$ of $M$ by leaf closures, but also valid for $\mfg\cdot \F$. In fact, the proof makes essential use of Theorem \ref{thm:Killing} below.

In Section \ref{sec:Molinobundle} we will recall the definition of a Killing foliation (see \cite{Mozgawa} or \cite[Section 5.5]{Molino}). Until then it is sufficient to know that a Riemannian foliation on a simply-connected manifold is automatically a Killing foliation. Isometric transverse actions appear naturally for Killing foliations due to the following theorem, which is a reformulation of  \cite[Theorem 5.2]{Molino} from our equivariant point of view. The details will be given in Section \ref{sec:Molinobundle}.
\begin{thm}\label{thm:Killing}
Let $\F$ be a Killing foliation on a complete manifold $M$. Then there exists an effective isometric transverse action of an abelian Lie algebra $\mfa$, called structural Killing algebra, such that $\mfa\cdot\F=\oF$.
\end{thm}
Note that from $\mfa\cdot\F=\oF$ it follows that the zero set $M^\mfa$ of the $\mfa$-action is the union of closed leaves. If the transverse action of $\mfa$ is free, then $\oF$ is a regular foliation. If it is transitive then all leaves are dense in $M$.

\section{Equivariant basic cohomology of transverse actions}\label{sec:equivbasiccoh}
The goal of this section is to introduce equivariant basic cohomology of a transverse action on a foliated manifold. As motivation, we first review the topological definition of equivariant cohomology of an action on a Lie group $G$ on a manifold via the Borel construction. In case $G$ is compact and connected, this object depends only on the induced infinitesimal action of the Lie algebra $\mfg$. To define the equivariant cohomology of an arbitrary infinitesimal action, the natural subsitute of the Borel construction is its well-known de-Rham-theoretic analogue in which the classifying space is replaced by the Weil algebra, see Example \ref{ex:weilmodel}. In case the action is induced by the action of a compact, connected Lie group, this so-called Weil model of equivariant cohomology coincides with the previous definition.

Following \cite{GS1999} we extend this definition of equivariant cohomology to $\mfg^\star$-algebras, as introduced by Cartan \cite[Section 4]{cartan2} (see also e.g.~\cite[Definition 2.3.1]{GS1999} or \cite[Definition 3.13]{kambertondeur}). We only introduce this definition as a motivation, as it will not be used in the subsequent sections. Instead, we will make constant use of yet another model, the Cartan model, which is isomorphic to the Weil model. 

Given a foliation $\F$ on a manifold $M$ with a transverse action of a Lie algebra $\mfg$, we will show in Proposition \ref{prop:omegaMFgstar} that the complex of basic differential forms $\Omega(M,\F)$ is a $\mfg^\star$-algebra. This allows us to define the equivariant basic cohomology $H^*_\mfg(M,\F)$ as a special case of the above general construction.

\subsection{Equivariant cohomology of Lie group actions: the Borel construction}
Let $M$ be a $G$-manifold. Let $EG$ be the classifying space of $G$, i.e., a contractible space on which $G$ acts freely. This space is unique up to $G$-homotopy equivalence. The projection $EG\to BG:=EG/G$ is called universal $G$-bundle. The Borel construction/homotopy quotient $M_G$ is defined as
$$
M_G:=EG\times_G M,
$$
and the equivariant cohomology $H^*_G(M)$ is defined as $H^*_G(M):=H^*(M_G)$, where the latter denotes singular cohomology with real coefficients.

\subsection{Equivariant cohomology of $\mfg^\star$-algebras} \label{sec:eqcohomgstar}

\begin{defn}\label{defn:g*}
Let $\mfg$ be a finite-dimensional Lie algebra and $A=\bigoplus A_k$ a $\ZZ$-graded algebra. We call $A$ a
(graded) $\mfg^\star$-algebra or differential graded $\mfg$-algebra if there is a derivation $d:A\to A$ of degree $1$ and  derivations $i_X:A\to A$ of degree $-1$ and $L_X:A\to A$ of degree $0$ for all $X\in \mfg$ such that:
\begin{enumerate}
\item $d^2=0$
\item $i_X^2=0$
\item $[L_X,L_Y]=L_{[X,Y]}$
\item $[L_X,i_Y]=i_{[X,Y]}$
\item $L_X=di_X+i_Xd$.
\end{enumerate}
\end{defn}

In \cite[Chapter 2]{GS1999} there is a reformulation of this definition phrased in the language of superalgebra:
Putting together the homomorphisms $i:\mfg\to {\rm Der}\, A; X\mapsto i_X$, $L:\mfg\to {\rm Der}\, A; X\mapsto L_X$ and the differential $d$, one obtains an action  of the Lie super algebra $\mfg^\star=\mfg\oplus\mfg\oplus\RR$, equipped with a natural Lie super bracket encoding the above identities, by superderivations on $A$.

\begin{ex} \label{ex:gstaraction}
An  infinitesimal action of a finite-dimensional  Lie algebra $\mfg$ on a manifold $M$ induces a $\mfg^\star$-algebra structure on the de Rham complex $\Omega^*(M)$.
\end{ex}

\begin{ex}
The Weil algebra $W(\mfg)=\Lambda(\mfg^*)\otimes S(\mfg^*)$ of a finite-dimensional  Lie algebra $\mfg$, the tensor product of the exterior and the symmetric algebra of $\mfg^*$, is an acyclic $\mfg^\star$-algebra, i.e., $H(W(\mfg),d)=\RR$. See \cite[Theorem 3.2.1]{GS1999}, \cite{cartan} or \cite[Lemma 4.3]{kambertondeur} for the definitions of $i_X,L_X,d$ and the proof.
\end{ex}

\begin{defn}
Let $A$ be a $\mfg^\star$-algebra. Then $A_{{\rm bas}\, \mfg}=\{\omega\in A\mid i_X\omega=L_X\omega=0\ \mbox{for all}\ X\in\mfg\}$ is $d$-invariant. It is called the {\em basic subcomplex} of $A$.
\end{defn}

We define the equivariant cohomology of an arbitrary $\mfg^\star$-algebra $A$ as
\begin{equation}\label{eq:weildef}
H_\mfg^*(A)=H((W(\mfg)\otimes A)_{{\rm bas\,}\mfg},d).
\end{equation}
Here the $\mfg^\star$-algebra structure on the tensor product is naturally extended from the factors.

\begin{ex}\label{ex:weilmodel}
Consider an infinitesimal action of a finite-dimensional  Lie algebra $\mfg$ on a manifold $M$ given by a homomorphism $\tau:\mfg\to \Xi(M)$. Let $W(\mfg)=\Lambda(\mfg^*)\otimes S(\mfg^*)$ be the Weil algebra of $\mfg$. Its basic complex $B\mfg=W(\mfg)_{{\rm bas\,}\mfg}$ is the ring of $\mfg$-invariant polynomials $S(\mfg^*)^\mfg$. In analogy to the Borel construction, the equivariant cohomology $H^*_\mfg(M)$ is defined as
$$
H^*_\mfg(M):=H^*_\mfg(\Omega(M))=H((W(\mfg)\otimes \Omega(M))_{{\rm bas\,}\mfg},d).
$$
If the action is induced by the action of a compact, connected Lie group $G$ we have
\begin{equation}\label{eq:inf}
H_G^*(M)=H^*_\mfg(M)
\end{equation}
by the equivariant de Rham Theorem, see e.g.~\cite[Theorem 2.5.1]{GS1999}. In this case $H_\mfg(M)$ is called the {\em Weil model} of equivariant cohomology.
\end{ex}

Another way to define the equivariant cohomology of an arbitrary $\mfg^\star$-algebra $A$ is the following. First define the {\em Cartan complex}
\begin{equation*}
C_\mfg(A):=(S(\mfg^*)\otimes A)^\mfg.
\end{equation*}
Here the superscript denotes the subspace of $\mfg$-invariant elements, i.e., those $\omega\in S(\mfg^*)\otimes A$ for which $L_X\omega=0$ for all $X\in \mfg$. The differential $d_\mfg$ of the {\em Cartan complex} $C_\mfg(A)$ is defined by
$$
(d_\mfg \omega)(X)=d(\omega(X))+i_X(\omega(X)),
$$
where we consider an element in $C_\mfg(A)$ as a $\mfg$-equivariant polynomial map $\mfg\to A$, see \cite[p.~45]{GS1999}.
Choosing a basis $\{X_i\}_{i=1,\ldots,r}$ of $\mfg$ with dual basis $\{u_i\}_{i=1,\ldots,r}$ of $\mfg^*$ we have
$$
d_\mfg\omega:=d\omega+\sum_{i=1}^ri_{X_i}(\omega)u_i.
$$
Now the equivariant cohomology of the $\mfg^\star$-algebra $A$ is defined as
\begin{equation}\label{eq:cartandef}
H_\mfg(A):=H(C_\mfg(A),d_\mfg).
\end{equation}
There is a natural $S(\mfg^*)^\mfg$-algebra structure on $H^*_\mfg(A)$, see \cite[Section 6.6]{GS1999}. This so-called {\em Cartan model} of equivariant cohomology is isomorphic to the former definition \eqref{eq:weildef} via the Mathai-Quillen isomorphism, see \cite[Section 4.1 and 4.9.2]{GS1999}.  From now on we will only use the Cartan model.

For the following definition compare \cite[Definition 2.3.4]{GS1999} and \cite[Appendix B, p.~182]{GLS1996}.
\begin{defn}
Let $A$ be a $\mfg^\star$-algebra. Let $X_i\in\mfg$ be a basis of $\mfg$. Then $A$ is  {\em free} if there are $\theta_i\in A_1$ such that $i_{X_j}(\theta_i)=\delta_{ij}$. If in addition the $\theta_i$ can be chosen such that their span in $A_1$ is $\mfg$-invariant, then $A$ is said to be {\em of type (C)}.
\end{defn}
Note that a free $\mfg^\star$-algebra $A$ is automatically of type (C) if the action of $\mfg$ on $A$ is induced by an action of a compact Lie group.
\begin{ex} If an infinitesimal isometric action of a finite-dimensional  Lie algebra $\mfg$ on a manifold $M$ is free, then the $\mfg^*$-algebra $\Omega(M)$ is of type (C).
\end{ex}

The following proposition generalizes the fact that the equivariant cohomology of a (locally) free Lie group action is the cohomology of the orbit space, see \cite[Appendix C.2]{GGK}.
\begin{prop}\label{prop:freebasic}
Let $A$ be a $\mfg^\star$-algebra of type (C). Then $H^*_\mfg(A)=H^*(A_{{\rm bas}\, \mfg})$.
\end{prop}
\begin{proof}
See \cite[Section 5.1]{GS1999} or \cite[Proposition B.8.1]{GLS1996}.
\end{proof}
There is a proof of the following proposition using the Weil model, see \cite[Section 4.6]{GS1999}. We think it is helpful to see how it can be proven in the Cartan model.
\begin{prop}\label{prop:abstractcommutingaction}
Let $A$ be an $(\mfh\times\mfk)^\star$-algebra with $A_k=0$ for $k<0$, which is of type (C) as an $\mfh^\star$-algebra. If either $A^\mfk=A$ or the $\mfk^\star$-algebra structure on $A$ is compatible with the action of a compact connected Lie group in the sense of \cite[Definition 2.3.1]{GS1999}, then
$$
H^*_{\mfh\times\mfk}(A)=H^*_{\mfk}(A_{{\rm bas\,}\mfh})
$$
as $S(\mfk^*)$-algebras.
\end{prop}
\begin{proof}
The two sides of the equation are the cohomologies of the complexes
$$
C=((S(\mfk^*)\otimes S(\mfh^*)\otimes A)^{\mfk\times \mfh},d_{\mfh\times \mfk})$$
and
$$
C'=((S(\mfk^*)\otimes A_{{\rm bas\,}\mfh})^\mfk,d_{\mfk}),
$$
respectively. We have a natural inclusion $i:C'\hookrightarrow C$ of complexes and thereby obtain a map $i_\sharp:H^*_{\mfk}(A_{{\rm bas\,}\mfh})\to H^*_{\mfh\times\mfk}(A)$.

We identify $C=(S(\mfk^*)\otimes (S(\mfh^*)\otimes A)^\mfh)^\mfk$ with the space of $\mfk$-invariant polynomials $\mfk\to (S(\mfh^*)\otimes A)^\mfh$ and define $\delta:C\to C$ via $\delta(\omega)(X)= i_X(\omega(X))$ for all $X\in \mfk$. We have $d_{\mfh\times \mfk}=d_\mfh+\delta$ and, restricted to $C'$, also $d_\mfk=d+\delta$. Introducing the bigradings $C^{p,q}=(S^p(\mfk^*)\otimes ((S(\mfh^*)\otimes A)^\mfh)^{q-p})^{\mfk}$ and $C'^{p,q}=(S^p(\mfk^*)\otimes A^{q-p}_{{\rm bas\,}\mfh})^\mfk$, we may regard $C$ and $C'$ as double complexes with vertical operators $d_\mfh$ and $d$ respectively, and $\delta$ as horizontal operator.

The spectral sequences of these double complexes satisfy
$$
E_1(C)=(S(\mfk^*)\otimes H^*_\mfh(A))^\mfk=(S(\mfk^*)\otimes H^*(A_{{\rm bas}\,\mfh}))^\mfk=E_1(C')
$$
by Proposition \ref{prop:freebasic} and the assumption on the $\mfk^\star$-algebra structure: if $A^\mfk=A$, then e.g.~$C=S(\mfk^*)^\mfk\otimes (S(\mfh^*)\otimes A)^\mfh$ and hence $E_1(C)= S(\mfk^*)^\mfk \otimes H^*_\mfh(A)$. In case the $\mfk^\star$-algebra structure is compatible with the action of a compact connected Lie group, we can use an averaging argument to calculate $E_1(C)$ and $E_1(C')$ as in \cite[Theorem 6.5.1]{GS1999}. The spectral sequences converge to $H^*(C)$ respectively $H^*(C')$ because of the assumption $A_k<0$ for $k<0$. Hence \cite[Theorem 6.4.2]{GS1999} the natural map $i_\sharp$ is an isomorphism.
\end{proof}
\subsection{Basic cohomology}
Let $(M,\F)$ be a foliated manifold of codimension $q$. A form $\omega\in\Omega(M)$ is {\em basic} if $i_X\omega=0$ and $L_X\omega=0$ for all $X\in\Xi(\F)$. We denote the space of basic $k$-forms by $\Omega^k(M,\F)$. It is a sub-$\Omega^0(M,\F)$-module of $\Omega(M)$. The space $\Omega^*(M,\F)$ equipped with the exterior product, is the {\em algebra of $\F$-basic differential forms}. It is $d$-invariant. Its cohomology, denoted by $H^*(M,\F)$ is the {\em basic cohomology algebra} of $(M,\F)$. We call $\dim H^k(M,\F)$ the $k$th-{\em basic Betti number}.

On a foliation chart domain $U$ with projection $p:U\to U/\F$ such that the foliation on $U$ is given by the fibers of $p$, the pullback $p^*:\Omega^*(U/\F)\to\Omega^*(U,\F)$ is an isomorphism of differential complexes. In particular $\Omega^k(M,\F)=0$ and therefore $H^k(M,\F)=0$ for $k>q$.
\begin{ex}
Let $(M,\F)$ be a simple foliation, i.e., there is a locally trivial fibration $p:M\to B$ of manifolds such that the leaves of $M$ are given as the fibers of $p$. Then $p^*:\Omega^*(B)\to\Omega^*(M,\F)$ is an isomorphism of differential complexes and $H^*(B)\cong H^*(M,\F)$.
\end{ex}
\begin{prop}\label{prop:basiccohomfinite} Let $\F$ be a Riemannian foliation on a complete manifold $M$ such that $M/\oF$ is compact. Then $\dim H^*(M,\F)<\infty$.
\end{prop}
\begin{proof} This result was proven in \cite[Th\'{e}or\`{e}me 0]{AlSeHe} under the more restrictive assumption of $M$ being compact. Under the assumption of $M/\oF$ being compact their proof is still valid as one may replace in \cite[Corollaire II.5]{AlSeHe} the notion of a transverse parallelism on a compact manifold by a complete transverse parallelism on a complete manifold (see Section \ref{sec:Molinobundle}) and because their space $W$ (which coincides with our space $W$ to be introduced in Section \ref{sec:Molinobundle}) is still compact due to $M/\oF=W/S$, see Theorem \ref{thm:equivariantMolino}.
\end{proof}
For more on this topic, see Chapter 2.3 and Appendix B by Sergiescu in \cite{Molino}.

\subsection{Equivariant basic cohomology}
The next observation is fundamental for this paper.
\begin{prop}\label{prop:omegaMFgstar}
A transverse action of a finite-dimensional  Lie algebra $\mfg$ on a foliated manifold $(M,\F)$ induces the structure of a $\mfg^\star$-algebra on $\Omega(M,\F)$.
\end{prop}
\begin{proof} Recall that a transverse action is a homomorphism $\mfg\to l(M,\F);X\mapsto X^*$. The basic complex $\Omega(M,\F)$ inherits the differential $d$ from $\Omega(M)$. For $X\in\mfg$ and $\omega\in\Omega(M,\F)$ let $i_X\omega=i_{\tilde X}\omega$ and $L_X\omega=L_{\tilde X}\omega$, where $\tilde X\in L(M,\F)$ is any foliate field whose transverse field is $X^*$. Note that these operators are well-defined because $\omega$ is basic. We are left to show that $i_{X}\omega$ and $L_X\omega$ are in $\Omega(M,\F)$ again, as the commutation relations in Definition \ref{defn:g*} are inherited from $\Omega(M)$. Let $A\in \Xi(\F)$. Then
$$
i_Ai_{X}\omega=i_Ai_{\tilde X}\omega=-i_{\tilde X}\underbrace{i_A\omega}_{=0}=0
$$
and
$$
L_Ai_{X}\omega=L_Ai_{\tilde X}\omega=i_{\tilde X}\underbrace{L_A\omega}_{=0}+i_{[A,\tilde X]}\omega=0.
$$
Here the last summand is zero since $[A,\tilde X]\in \Xi(\F)$. So $i_{X}\omega\in  \Omega(M,\F)$. Similarly $L_{X}\omega\in \Omega(M,\F)$.
\end{proof}
This proposition enables us to apply the general construction of equivariant cohomology of a $\mfg^*$-algebra (via the Cartan model), see \eqref{eq:cartandef}, to $\Omega(M,\F)$. \begin{defn} For a foliated manifold $(M,\F)$ with a transverse action of a finite-dimensional  Lie algebra $\mfg$ we write
$$
\Omega_\mfg(M,\F):=C_\mfg(\Omega(M,F))=(S(\mfg^*)\otimes \Omega(M,\F))^\mfg
$$
for the Cartan complex associated to the $\mfg^*$-algebra $\Omega(M,\F)$ and call it the complex of ($\mfg$-)equivariant ($\F$-)basic forms. The {\em equivariant basic cohomology} of the $\mfg$-action on $(M,\F)$ is
 $$
H_\mfg(M,\F):=H_\mfg(\Omega(M,\F))=H(\Omega_\mfg(M,\F),d_\mfg).
$$
\end{defn}

\begin{ex}
Let $\pi:M\to B$ a locally trivial fibration with connected fibers. We consider the set of fibers as a foliation $\F$ on $M$. Assume that there is an infinitesimal action on $B$ by a finite-dimensional  Lie algebra $\mfg$. Then we can define a transverse action on $(M,\F)$ by $\mfg$ such that $\pi$ becomes $\mfg$-equivariant with respect to these (transverse) actions. Then
$$
H_\mfg(B)=H_\mfg(M,\F).
$$
If moreover the action of $\mfg$ on $B$ is integrable to an action of a compact connected Lie group $G$ then $H_{G}(B)=H_\mfg(B)=H_\mfg(M,\F).$
\end{ex}

\begin{lemma} \label{lem:basicformsainvariant}Let $\F$ be a Killing foliation on a complete manifold $M$ with structural Killing algebra $\mfa$. Then for any $\F$-saturated open (hence $\mfa$-invariant) subset $U\subset M$, we have $\Omega(U,\F)^\mfa=\Omega(U,\F)$. In particular, $\Omega_\mfa(U,\F)=\Omega(U,\F)\otimes S(\mfa^*)$.
\end{lemma}
\begin{proof} We use the notation of \cite[Appendix D by E.~Salem]{Molino}. See also \cite{Salem1988}. Let $\{U_i\}$ be a covering of $M$ by foliation chart domains with transversals $T_i\subset U_i$. Then $\mfa$ acts isometrically on the Riemannian manifolds $T_i$. For small elements $X\in \mfa$, the transformation $\exp X$ is an element of the closure of the holonomy pseudogroup of $\F$, see \cite[p.~287]{Molino}. Thus the basic forms are $\mfa$-invariant.
\end{proof}

\begin{prop}\label{prop:mayervietoris} Let $\F$ be a Killing foliation on a complete manifold $M$ with structural Killing algebra $\mfa$.
Let $U$ and $V$ be $\F$-saturated open subsets of $M$. Then the Mayer-Vietoris sequence
$$
\ldots \longrightarrow H^*_\mfa(U\cup V,\F)\longrightarrow H^*_\mfa(U,\F)\oplus H^*_\mfa(V,\F) \longrightarrow H^*_\mfa(U\cap V,\F)\longrightarrow \ldots
$$
is exact.
\end{prop}
\begin{proof}
We have that
$$
0\longrightarrow\Omega^*(U\cup V,\F)\longrightarrow \Omega^*(U,\F)\oplus \Omega^*(V,\F) \longrightarrow \Omega^*(U\cap V,\F)\longrightarrow 0
$$
is a short exact sequence of $\mfa^\star$-modules. By Lemma \ref{lem:basicformsainvariant}, tensoring with $S(\mfa^*)$ gives the short exact sequence 
$$
0\longrightarrow\Omega_\mfa^*(U\cup V,\F)\longrightarrow \Omega_\mfa^*(U,\F)\oplus \Omega_\mfa^*(V,\F) \longrightarrow \Omega_\mfa^*(U\cap V,\F)\longrightarrow 0
$$
which induces the long exact Mayer-Vietoris sequence. Note that the argument in \cite[Section 4.8]{GS1999} cannot be used for transverse actions of arbitrary Lie algebras as for infinitesimal actions there is no averaging process available.
\end{proof}

\subsection{Basic calculations} \label{sec:basiccalculations}

\begin{lemma}\label{lem:closedleaves}
Let $\F$ be a Riemannian foliation on a complete manifold $M$ with transverse metric $g$ and $U\subset M$ be $\F$-saturated and open. If all leaves of $\F$ in $U$ are closed then $H^*(U,\F)=H^*(U/\F)$.
\end{lemma}
\begin{proof}
Let $\pi:M\to M/F$ be the projection. By \cite[Proposition 3.7]{Molino} $M/\F$ carries the structure of an orbifold   (note that the assumption of compact leaves can be relaxed to require only closed leaves).  The map $\pi^*:\Omega(U/\F)\to\Omega(U,\F)$ is an isomorphism of differential complexes, where the left object is the space of orbifold differential forms. $H^*(U,\F)=H^*(U/\F)$ follows. The right hand side may be interpreted as the cohomology of $\Omega(U/\F)$ or the singular cohomology of $U/\F$ by \cite[Theorem 1]{Satake1956}.
\end{proof}
\begin{lemma}\label{lem:freetransverseaction}
Let $\F$ be a Riemannian foliation on a manifold $M$ with transverse metric $g$. Consider an isometric transverse action of a Lie algebra $\mfg$ on $(M,\F)$. If the $\mfg$-action is free then the $\mfg^\star$-algebra $\Omega(M,\F)$ is of type (C) and $H_\mfg(M,\F)=H^*(M,\mfg \F)$. If $M$ is complete and all leaves of $\mfg\F$ in $U$ are closed then for every $\F$-saturated, open subset $U\subset M$ we have $H_\mfg(U,\F)=H(U/\mfg\F)$.
\end{lemma}
\begin{proof}
Recall from the proof of Proposition \ref{prop:omegaMFgstar} the definition of the $\mfg^\star$-algebra structure on $\Omega(M,\F)$: for $X\in \mfg$ we have $i_X\omega=i_{\tilde X}\omega$ and $L_X\omega=L_{\tilde X}\omega$, where $\tilde X\in L(M,\F)$ is any foliate field whose transverse field is $X^*$. Let $X_i$ be a basis of $\mfg$ and let $\theta_i\in\Omega^1(M)$ such that $i_{\tilde X_i}\theta_j=\delta_{ij}$ and $\ker \theta_i(p)\supset T_p\F\oplus (\mfg\cdot p)^\perp$. Here, $(\mfg\cdot p)^\perp$ is the $g$-orthogonal complement of $\mfg\cdot p$ in $T_pM/T_p\F$, and thus $T_p\F\oplus (\mfg\cdot p)^\perp$ is in a natural way a subspace of $T_pM$. We claim $\theta_i\in \Omega^1(M,\F)$. By definition $\theta_i$ is horizontal. For $X\in\Xi(\F)$ we have
\begin{equation*}
(L_X\theta_i)(Y)=L_X\theta_i(Y)-\theta_i([X,Y])
\end{equation*}
for all vector fields $Y$. For $Y=\tilde X_j$ the first summand is  $X\delta_{ij}=0$ and the second vanishes because $[X,\tilde X_j]\in\Xi(\F)$ as $\tilde X_j\in L(M,\F)$. Now let $Y$ be transversely orthogonal to the $\mfg$-orbits, i.e., $Y_p\in T_p\F\oplus (\mfg\cdot p)^\perp$ for all $p$. Then the first summand vanishes and for the second observe that because of
\begin{align*}
0&=(L_Xg)(\tilde X_j,Y)\\
&=Xg(\tilde X_j,Y)-g([X,\tilde X_j],Y)-g(\tilde X_j,[X,Y])=-g(\tilde X_j,[X,Y]),
\end{align*}
$[X,Y]$ is transversely orthogonal to the $\mfg$-orbits. In this calculation we considered $g$ as a degenerate metric on $TM$ with kernel $T\F$. This proves the claim.

Next we will show that  the span of the $\theta_i$, the connection elements, is $\mfg$-invariant. We have
$$
(L_{X_j}\theta_i)(Y)=(L_{\tilde X_j}\theta_i)(Y)=L_{\tilde X_j}\theta_i(Y)-\theta_i([\tilde X_j,Y])
$$
Both summands on the right hand side are zero for $Y\in\Xi(\F)$ because $\tilde X_j$ is foliate, and for $Y$ transversely orthogonal to the $\mfg$-orbits because the space of transversely orthogonal vector fields is $\mfg$-invariant by a calculation using $L_{\tilde X_k}g=0$ as before. So $\ker L_{X_j}\theta_i(p)\supset T_p\F\oplus (\mfg\cdot p)^\perp$. For $Y=\tilde X_k$ the first summand is zero and the second is constant because $\mfg$ is a Lie algebra. Thus $L_{X_j}\theta_i$ is a linear combination of the $\theta_j$. It follows that $\Omega(M,\F)$ is of type (C).

Then by Proposition \ref{prop:freebasic} we have $H^*_\mfg(\Omega(M,\F))=H^*(\Omega(M,\F)_{{\text{bas\,}\mfg}})$.  The fact that $\Omega(M,\F)_{{\text{bas\,}\mfg}}$ consists of those elements in $\Omega(M)$ that are horizontal and invariant with respect to $\Xi(\F)$ and $\mfg$ implies $\Omega(M,\F)_{{\text{bas\,}\mfg}}=\Omega(M,\mfg\F)$. Now the last statement of the Lemma follows from Lemma \ref{lem:closedleaves} because if $\mfg$ acts freely then $\mfg\cdot \F$ is a Riemannian foliation.
\end{proof}
Now let $\F$ be a Killing foliation on a complete manifold $M$ with structural algebra $\mfa$. Recall from Theorem \ref{thm:Killing} that $\mfa$ is abelian and $\mfa\cdot \F=\oF$. Together with the previous result  we obtain:
\begin{cor}\label{cor:afree}
If the structural algebra $\mfa$ of $(M,\F)$ acts freely on an open, $\F$-saturated subset $U\subset M$, then $H_\mfa^*(U,\F)=H^*(U/\oF)$.
\end{cor}
The next result is a generalization.
\begin{cor}\label{cor:oneisotropy}
If $\mfa$ acts with only one isotropy algebra $\mfk$ on an open, $\F$-saturated subset $U\subset M$, then $H_\mfa(U,\F)=H^*(U/\oF)\otimes S(\mfk^*)$.
\end{cor}
\begin{proof}
Let $\mfk'$ be a complement of $\mfk$ in $\mfa$.
We have $\Omega_\mfa(U,\F)=S(\mfa^*)\otimes \Omega(U,\F)^\mfa=S(\mfa^*)\otimes \Omega(U,\F)=S(\mfk^*)\otimes (S(\mfk'^*)\otimes \Omega(U,\F))$ by Lemma \ref{lem:basicformsainvariant}. 
The equivariant differential $d_\mfa$ leaves the space $S(\mfk'^*)\otimes \Omega(U,\F)$ invariant and is zero on $S(\mfk^*)$. Thus, $
H^*_\mfa(U,\F)=S(\mfk^*)\otimes H^*(S(\mfk'^*)\otimes \Omega(U,\F),d_{\mfk'})=S(\mfk^*)\otimes H^*_{\mfk'}(U,\F)=S(\mfk^*)\otimes H^*(U,\mfk'\F)$; the last equality is due to Lemma \ref{lem:freetransverseaction}. Since $\mfk'\F=\mfa\F=\oF$ the claim follows.
\end{proof}
For any leaf $L$ we may apply the corollary to $\F$, restricted to the complete manifold $\oL$.
\begin{cor}\label{cor:cohomologyleafclos}
For every leaf $L$ of $\F$ we have $H^*_\mfa({\overline L},\F)=S(\mfa_L^*)$.
\end{cor}
Here, $\mfa_L$ is the unique $\mfa$-isotropy algebra of the $\mfa$-action on $L$, see Section \ref{sec:stratification}.

\subsection{Equivariantly formal actions} The results from this subsection are transferred from \cite[Section 6.5]{GS1999} to our case. Let $\mfg$ be a transverse action on a foliated manifold $(M,\F)$.

If the transverse action of $\mfg$ is trivial, then we have $H_\mfg^*(M,\F)= S(\mfg^*)^\mfg\otimes H^*(M,\F)$ as $S(\mfg^*)^\mfg$-modules (even as $S(\mfg^*)^\mfg$-algebras). In analogy to the case of Lie group actions on manifolds we define
\begin{defn}
If $H_\mfg^*(M,\F)\cong S(\mfg^*)^\mfg\otimes H^*(M,\F)$ as $S(\mfg^*)^\mfg$-modules, we call the $\mfg$-action {\em equivariantly formal}.
\end{defn}
For Lie group actions, this property was studied in detail in \cite{GKM}. See also \cite[Section 6.9]{GS1999}. From now on, let $(M,\F)$ be a Killing foliation with structural Killing algebra $\mfa$. Recall  $\Omega_\mfa(M,\F)=S(\mfa^*)\otimes \Omega(M,\F)$ by  Lemma \ref{lem:basicformsainvariant}.

\begin{thm}\label{thm:specseq}
The spectral sequence of the double complex
$$
\Omega^{p,q}_\mfa(M,\F):=S^p(\mfa^*)\otimes \Omega^{q-p}(M,\F)
$$
with vertical operator $d$ and horizontal operator $\delta$ defined by
$\delta(\omega)(X)=i_X(\omega(X))$
has
$$
E_1^{p,q}=S^p(\mfa^*)\otimes H^{q-p}(M,\F)
$$
and converges to $H^*_\mfa(M,\F)$.
\end{thm}
 Note that for the convergence we require the fact $H^d(M,\F)=0$ for all $d>\codim\F$. For a general transverse action of a Lie algebra $\mfg$ the spectral sequence of the double complex $\Omega^{p,q}_\mfg(M,\F)=(S^p(\mfg^*)\otimes \Omega^{q-p}(M,\F))^\mfg$ still exists, but the shape of the $E_1$-term is unclear as it might happen that $\Omega(M,\F)^\mfg\subsetneq \Omega(M,\F)$ and there is no averaging process available (which is needed to imitate the proof of \cite[Theorem 6.5.1]{GS1999}). 
 
As in \cite[Theorem 6.6.1]{GS1999}, Theorem \ref{thm:specseq} implies
\begin{prop}\label{prop:finitelygenerated}
If $\dim H^*(M,\F)<\infty$, then $H^*_\mfa(M,\F)$ is a finitely generated $S(\mfa^*)$-module.
\end{prop}
Combining this with Proposition \ref{prop:basiccohomfinite}, we obtain
\begin{cor}\label{cor:leafspacecptthenfinite}
If $M/\oF$ is compact, then $H^*_\mfa(M,\F)$ is a finitely generated $S(\mfa^*)$-module.
\end{cor}
\begin{prop}\label{prop:E2collapse}
The transverse action is equivariantly formal if and only if the spectral sequence of the double complex $\Omega^{p,q}_\mfa(M,\F)$ collapses at the $E_1$-stage.
\end{prop}
\begin{proof}
This is a standard argument based on the rank of the $E_k$-terms of the spectral sequence as modules over $S(\mfa^*)$, see e.g.~\cite[Lemma C.24]{GGK}.
\end{proof}
\begin{cor}[{\cite[Theorem 6.5.3]{GS1999}}]\label{cor:odd}
If $H^{\rm odd}(M,\F)=0$ then the transverse action of $\mfa$ is equivariantly formal.
\end{cor}
\begin{cor}
The transverse action is equivariantly formal if and only if the natural map $H^*_\mfa(M,\F)\to H^*(M,\F)$ is surjective.
\end{cor}
\begin{proof}
By Proposition \ref{prop:E2collapse}, equivariant formality of the transverse action is equivalent to the vanishing of all higher differentials $d_1=\delta,d_2,d_3,\ldots$ in the spectral sequence of the double complex $\Omega^{p,q}_\mfa(M,\F)$. By definition of the higher differentials (see e.g.~Section 6.1 in \cite{GS1999}) it is thus clear that in case the action is equivariantly formal, every closed form in $H^*(M,\F)$ can be extended to a closed equivariant form.

If conversely $H^*_\mfa(M,\F)\to H^*(M,\F)$ is surjective, then $d_1$ vanishes on $E_1^{0,*}$. As $d_1$ vanishes on $S(\mfg^*)^\mfg=E_1^{*,0}$ as well, it is linear with respect to the $S(\mfa^*)$-module structure of $E_1$, hence $d_1=0$. Thus $E_2=E_1$ and the same argument applies inductively to the higher differentials. See also the proof of \cite[Theorem 5.10]{McC}.
\end{proof}

For standard torus actions, the following criterion for equivariant formality was proven in \cite[Corollary (4.2.3)]{allday}.

\begin{prop} \label{prop:freemodequivformal}
The transverse action of $\mfa$ is equivariantly formal if and only if $H^*_\mfa(M,\F)$ is a free $S(\mfa^*)$-module.
\end{prop}
\begin{proof}
By definition of equivariant formality, we need to show that if $H^*_\mfa(M,\F)$ is a free $S(\mfa^*)$-module, then $\dim H^*_\mfa(M,\F)\otimes_{S(\mfa^
*)}\RR =\dim H^*(M,\F)$.

The proof carries over almost verbatim from \cite[Proposition (3.3.15)]{allday}. We prove that without any assumption on $H^*_\mfa(M,\F)$ there is an Eilenberg-Moore spectral sequence converging to $H^*(M,\F)$ with $E_2$-term equal to ${\rm Tor}_*^{S(\mfa^*)}(H^*_\mfa(M,\F),\RR)$. Compare \cite[Example 1, p.~38]{Hsiang} for the Eilenberg-Moore spectral sequence in the case of a Lie group action.

Consider the fourth quadrant double complex
\[
C^{p,q}=\Omega_\mfa^{p+2q}(M,\F)\otimes \Lambda^{-q}(\mfa^*),
\]
where $p\geq 0$ and $q\leq 0$, with differentials $d_1:C^{p,q}\to C^{p+1,q}$ given by $d_1(\omega\otimes 1)=d_\mfa\omega\otimes 1$ and $d_1(1\otimes s)=0$ for all $s\in \mfa^*$, and $d_2:C^{p,q}\to C^{p,q+1}$ given by $d_2(1\otimes s)=s\otimes 1$ and $d_2(\omega\otimes 1)=0$. We obtain two spectral sequences with $E_2$-terms $H(H(C^{*,*},d_1),d_2)$ and $H(H(C^{*,*},d_2),d_1)$, respectively, both converging to $H(C^{*,*},d_1+d_2)$. Let us first consider the latter.

Recall that if $r=\dim \mfa$, then
\[
0\to S(\mfa^*)\otimes \Lambda^r(\mfa^*)\to\ldots\to S(\mfa^*) \otimes \Lambda^1(\mfa^*)\to S(\mfa^*)\to \RR\to 0,
\]
where the maps are defined analogously to $d_2$, is a free resolution of the $S(\mfa^*)$-module $\RR$, called Koszul resolution.
Tensoring this sequence with $\Omega^*_\mfa(M,\F)$ and deleting the last term, we obtain the complex
\[
0\to \Omega^{*-2r}_\mfa(M,\F)\otimes \Lambda^r(\mfa^*)\overset{d_2}{\to}\ldots\overset{d_2}{\to} \Omega^{*-2}_\mfa(M,\F) \otimes \Lambda^1(\mfa^*)\overset{d_2}{\to} \Omega^*_\mfa(M,\F)\to 0,
\]
whose cohomology on the one hand equals $H(C^{*,*},d_2)$, and on the other hand computes ${\rm Tor}_*^{S(\mfa^*)}(\Omega^*_\mfa(M,\F),\RR)$. As $\Omega_\mfa^*(M,\F)$ is a free $S(\mfa^*)$-module, $H(C^{*,*},d_2)$ consists of only one row, namely $\Omega_\mfa^*(M,\F)\otimes_{S(\mfa^*)}\RR$, whose $d_1$-cohomology equals $H^*(M,\F)$.

An analogous argument shows that the other spectral sequence has $E_2$-term equal to ${\rm Tor}_*^{S(\mfa^*)}(H^*_\mfa(M,\F),\RR)$. Because of freeness of $H^*_\mfa(M,\F)$ this is the same as $H^*_\mfa(M,\F)\otimes_{S(\mfa^*)}\RR$ and the spectral sequence therefore collapses at $E_2$. It follows $\dim H^*_\mfa(M,\F)\otimes_{S(\mfa^*)}\RR=\dim H^*(M,\F)$.
\end{proof}

\section{The correspondence principle}\label{sec:correspondence}
In this section we will review Molino's construction of an $\Oq$-manifold $W$ associated to a Riemannian foliation $\F$ on a complete manifold $M$ and find a correspondence between basic data of $(M,\F)$ and equivariant data of $(W,\Oq)$.
\subsection{The Molino bundle}\label{sec:Molinobundle}
Let $M$ be a complete, connected $n$-dimensional manifold with a Riemannian foliation $\F$ and transverse metric $g$. We denote the codimension of $\F$ by $q$. Let $P=\whM$ be the principal $S=\Oq$-bundle over $M$ of orthonormal frames of $\nu\F=TM/T\F$, the {\em Molino bundle} of $(M,\F)$. We denote the natural projection $P\to M$ by $\pi$. If $(M,\F)$ is transversely orientable then $P$ has two connected $\SOq$-invariant components. A choice of transverse orientation corresponds to a choice of a component. In this case we will by abuse of notation denote this component also by $P$ and let $S=\SOq$.  The normal bundle $\nu\F$ is associated to $P$, and we denote the connection form on $P$ corresponding to the basic Riemannian connection $\nabla$ on $\nu\F$ introduced in Section \ref{subsec:strkillalg} by $\omega$. We write $H_1:=\ker \omega$ for the $S$-invariant horizontal distribution. The manifold $P$ carries a natural foliation $\E$ obtained by horizontally lifting the leaves of $M$. This foliation is $S$-invariant, i.e., $S$ maps leaves to leaves.  By construction $\omega$ is a basic form with respect to $\E$, i.e., $i_X\omega=0$ and $L_X\omega=0$ for all $X\in \Xi(\E)$. In particular, we may regard $\omega$ as a map $\omega:\nu \E\to \soq$. We write $\cH:=H_1/T\E$ for the {\it transverse horizontal distribution}. Now we lift the transverse metric $g$ on $\nu\F$ to an $S$-invariant metric on the $S$-invariant distribution $\cH$ of $\nu\E$. Recall \cite[p.~70, p.~148]{Molino} the fundamental $1$-form $\theta:\nu\E\to \RR^q$ defined by 
$$
\theta(X_{\hat x})=\hat x^{-1}(\pi_*(X_{\hat x}))
$$
where $\hat x$ is an orthonormal frame of $\nu_x\F$, understood as the isomorphism $\hat x:\RR^q\to\nu_x\F$ sending the canonical basis to the frame, and $X_{\hat x}\in\nu_{\hat x}\E= T_{\hat x}P/T_{\hat x} \E$. The fundamental $1$-form $\theta$ is $\E$-basic by \cite[Lemma 2.1 (i)]{Molino}. By definition, the standard scalar product on $\RR^q$, pulled back with $\theta$ to a degenerate metric on $\nu \E$, coincides with $\pi^*g$.  We consider the $\E$-basic, $S$-equivariant map $\theta\oplus\omega:\nu\E\to\mR^q\oplus \soq$. Pulling back the sum of the standard scalar product on $\RR^q$ and an arbitrary (unique up to a scalar)  biinvariant metric on $\soq$ with $\theta\oplus \omega$ yields an  $S$-invariant $\E$-transverse metric $g_P$ on $(P,\E)$ with respect to which $\E$ is a Riemannian foliation. The projection $\pi$ becomes a transversely Riemannian submersion, i.e., $\pi_*:\nu\E\to\nu\F$ is surjective and the restriction $\pi_*:\cH\to\nu\F$ is orthogonal at any point.

Molino argued that the foliation $\E$ has a global transverse parallelism, i.e., $\nu\E$ is parallelizable by transverse fields (we say $\E$ is TP), see \cite[p.~82, p.~148]{Molino}. Using the metric $g_P$ constructed above we can slightly extend this argument to obtain a complete transverse parallelism (we say $\E$ is TC), i.e., the transverse fields giving the parallelization have complete representatives in $L(P,\E)$: We choose an orthonormal basis of $\mR^q\oplus \soq$ compatible with the direct sum decomposition. The pull-back gives an orthonormal transverse parallelism of $\E$. Moreover, since $P$ is complete (as $M$ is complete) and since these transverse fields have constant length by construction of the metric, they have complete representatives in $L(P,\E)$. 

Let $C(P,\E)$ be the sheaf of local transverse fields that commute with all global transverse fields, called {\em commuting sheaf} of the TP-foliation $\E$, see \cite[Section 4.4]{Molino}. Its stalk is a Lie algebra $\mfg$ and we write $\umfg:=C(P,\E)$. Note that by  \cite[Lemma 4.6]{Molino} a section $X\in\umfg(U)$ on a connected subset $U$ is determined by its germ $X_p\in \umfg_p\cong\mfg$ at a point $p$. As $\umfg(U)\subset l(U,\E)$ for $U\subset P$ we may say that $\umfg$ defines a local transverse action on $(P,\E)$. This transverse action is isometric as it commutes with any transverse field and in particular with the orthonormal transverse parallelism above.  Theorem 4.3 of \cite{Molino} says that its orbits are the leaf closures of $\E$, i.e.,
\begin{equation}\label{eq:orbitsareclosures}
\{X_p\mid X\in \umfg_p\}\oplus T_p\E=T_p\oE
\end{equation} for all $p$. For short, we write $\umfg\cdot\E=\oE$.   Any local section of $\umfg$ is the natural lift of a local transverse Killing field of $M$, see  \cite[Proposition 3.4]{Molino}. The {\em Molino sheaf} of $(M,\F)$ is defined as the sheaf on $M$ whose sections are the local transverse Killing fields that naturally lift to local sections of $\umfg$. Thus the Molino sheaf of $(M,\F)$ can be identified with the push-forward $\pi_*C(P,\E)$ of the commuting sheaf. By abuse of notation we also denote the Molino sheaf by $\umfg$. Its stalk is also $\mfg$ which we call the {\em structural Killing algebra} of $(M,\F)$. In Molino's terminology the Molino sheaf is called commuting sheaf of $(M,\F)$, denoted by $C(M,\F)$. A slight difference is that for him, the structural algebra of $(M,\F)$ is the inverse Lie algebra of the stalk of $C(M,\F)$.

Regarding $\umfg$ as a local isometric transverse action on $(P,\E)$ and $(M,\F)$, we see that by construction $\pi:P\to M$ is $\umfg$-equivariant. By Lemma 5.1 of \cite{Molino} $\pi$ maps leaf closures of $\E$ onto leaf closures of $\F$; we can change the argument given there based on the compactness of $M$ by using the properness of $\pi$. Then $\umfg\cdot\F=\oF$ follows.

Since $\E$ is TP the foliation $\oE$ by leaf closures is simple, i.e., $W:=P/\oE$ is a manifold and $\oE$ is given as the set of fibers of the locally trivial fibration $\rho:P\to W$ (\cite[Proposition 4.1']{Molino}), called the {\em basic fibration}. As the right action of $S$ respects $\E$ it also respects $\oE$ by continuity and therefore descends to an action on $W$, so that $\rho$ is $S$-equivariant. Since the transverse metric $g_P$ on $\nu\E$ is $S$- and $\umfg$-invariant by \cite[Proposition 3.4.(ii)]{Molino}, $W$ can be equipped with a metric $g_W$ such that $\rho$ becomes a transversely Riemannian submersion that is equivariant with respect to the isometric $S$-actions.

For the next definition note that the Molino sheaf is locally constant by \cite[Theorem 5.2]{Molino}.
\begin{defn}[\cite{Mozgawa}]\label{defn:Killing}
A {\em Killing foliation} is a Riemannian foliation whose Molino sheaf is globally constant.
\end{defn}
A Riemannian foliation on a simply-connected manifold is therefore automatically a Killing foliation. This important special case was investigated in \cite{Ghys}. If $(M,\F)$ is a Killing foliation, then $C(P,\E)$ is also globally constant and thus the structural Killing algebra $\mfg$ is the center of $l(P,\E)$. Identifying $\mfg$ with the Lie algebra of global sections of the Molino sheaf, $\mfg$ is contained and central in $l(M,\F)$. See also \cite[Theorem 5.2]{Molino}. Therefore $\mfg$ is abelian (but not necessarily the full center of $l(M,\F)$). In order to indicate this we will denote the structural Killing algebra of a Killing foliation by $\mfa$. Note that in the case of a Killing foliation the identities $\umfg\cdot \E=\oE$ and $\umfg\cdot \F=\oF$ from above (see \eqref{eq:orbitsareclosures}) are equivalent to  $\mfa\cdot \E=\oE$ and $\mfa\cdot \F=\oF$  in the sense of transverse actions from Section \ref{sec:transverseaction}. 
We sum up Molino's results for Killing foliations:
\begin{thm}\label{thm:equivariantMolino}
Let $\F$ be a Killing foliation on a complete connected manifold $M$ with structural Killing algebra $\mfa$. Then the isometric transverse action of $\mfa$ lifts to a free isometric transverse action of $(P,\E)$ which commutes with the right $S$-action. Its respective orbits are the leaf closures, i.e., $\mfa\cdot \E=\oE$ and $\mfa\cdot \F=\oF$. The projection $\pi:P\to M$ maps leaves to leaves, leaf closures to leaf closures and is $\mfa$-equivariant. The projection $\rho:P\to W$ is equivariant with respect to the isometric $S$-actions and we obtain the following commutative diagram:
$$
\xymatrix{
  (P,\E,\mfa,S)\ar[d]^\pi  \ar[r]^\varrho& (W,S)\ar[d] \\
  (M,\F,\mfa) \ar[r]& M/\oF=W/S.}
$$
\end{thm}
The vertical arrows mean "mod $S$", the horizontal "mod leaf closures" or "mod $\mfa$". The bijection $M/\oF=W/S$ will be explained in Section \ref{sec:stratification}.
\begin{ex}\label{ex:torusaction}
Consider an almost effective isometric action without fixed points of a torus $T$ on a Riemannian manifold $M$, and a connected, dense Lie subgroup $B\subset T$ that acts locally freely. Denote by $\F$ the Riemannian foliation given by the $B$-orbits. Note that  $M/\oF=M/T$. The lifted action of $T$ on $P$ is free.  A $B$-equivariant (and hence $T$-equivariant, because $B$ is dense in $T$) field on $TP$ commutes with the fields induced by the $\mfb$-action and thus is foliate. Conversely, a transverse field in $l(P,\E)$ is $B$-equivariant because it can be represented by a field with values in the $T$-invariant geometric normal bundle of $\E$ defined by the natural Riemannian metric on $P$. Thus, $l(P,\E)$ is the same as the space of $T$-equivariant fields on $TP/T\E$.

It follows that the natural map $\mft/\mfb\to l(P,\E)$ maps $\mft/\mfb$ into the center of $l(P,\E)$. The dimension of the structural Killing algebra is by \eqref{eq:orbitsareclosures} the difference of dimensions of a leaf and its closure. For any $p\in P$ the closure of the leaf through $p$ is $Tp$, and thus this difference is $\dim Tp-\dim\E=\dim\mft/\mfb$. It follows that $\F$ is a Killing foliation and $\mft/\mfb$ is isomorphic to the structural Killing algeba $\mfa$. 

We furthermore have
\begin{equation*}
\begin{aligned}
H^*_T(M)&=H^*_{\mfa\oplus\mfb}(M)\\
&=H_\mfa^*(\Omega(M)_{{\rm bas\,}\mfb})\\
&=H_\mfa^*(\Omega(M,\F))&\\
&=H_\mfa^*(M,\F).&
\end{aligned}
\qquad
\begin{aligned}
&\mbox{as}\ T\ \mbox{is compact}\\
&\mbox{by Proposition}\ \ref{prop:abstractcommutingaction}\\ \, \\ \,
\end{aligned}
\end{equation*}
\end{ex}

\begin{ex} \label{ex:hopfexample} 
As a concrete example for the general construction in Example \ref{ex:torusaction} we consider the action of the two-dimensional torus $T^2$ on $S^3=\{(z_1,z_2)\in \CC^2\mid |z_1|^2+|z_2|^2=1\}$ given by
\[
((t_1,t_2),(z_1,z_2))=(t_1z_1,t_2z_2).
\]
For rationally independent reals $s_1,s_2$, the image $B$ of the inclusion $\RR\to T^2;\, t\mapsto (e^{2\pi i s_1t},e^{2\pi i s_2 t})$ is dense in $T^2$ and acts locally freely. We refer to the foliation of $S^3$ by $B$-orbits as the Hopf flow. It has exactly two closed orbits. In later sections we will illustrate our methods with this example.
\end{ex}

\subsection{The $\mfa$-stratification of $(M,\F)$}\label{sec:stratification}
We will describe the infinitesimal isotropy stratifcation of the Killing foliation $(M,\F)$. Let $L$ be a leaf. The isotropy algebra $\mfa_x=\{X\in\mfa\mid X^*_x=0\}$ is the same for all points $x\in L$; we set $\mfa_L:=\mfa_x$. In fact it is the same for all $x\in \oL=\mfa\cdot L$ and therefore also denoted by $\mfa_{\oL}$; this can be seen by passing to the natural $\mfa$-action on a local quotient $U/\F$ of a foliation chart domain $U$. For any subspace $\mfh\subset \mfa$, let $M^\mfh=\{p\in M\mid X^*_p=0 \text{ for all } X\in \mfh\}$ be the zero set of the transverse action of $\mfh$. For example, $M=M^{\{0\}}$ and $M^\mfa$ is the union $C$ of the closed leaves. Moreover we define $M^{(\mfh)}=\{p\in M\mid \mfa_p=\mfh\}\subset M^\mfh$. The {\em $\mfa$-stratification of $(M,\F)$} is the partition by the connected components of $M^{(\mfh)}$, where $\mfh$ runs over all isotropy algebras. An element of a stratification is called a {\em stratum}. For an arbitrary stratification the {\em  bottom stratum} is defined as the union of the closed strata. Each component of the union $C$ of the closed leaves is a component of the bottom stratum of the $\mfa$-stratification of $(M,\F)$ .

We will now define a stratification on the $S$-manifold $W$.   Given    a subalgebra $\mfh \subset \mfs$, we denote by $(\mfh)=\{\mfk\mid \Ad_s\mfh=\mfk\ \text{for some}\ s\in S\}$ the $S$-conjugacy class of $\mfh$ in $\mfs$. While for an orbit $Sx$ the isotropy group $\mfs_y$ depends on the choice of $y\in Sx$, the conjugacy class $(\mfs_x)$ does not and is therefore an invariant of $Sx$. The class $(\mfs_x)$ is called the {\em (infinitesimal) orbit type} of $Sx$. For an isotropy algebra $\mfh$, define the $(\mfh)$-{\it (infinitesimal) orbit type submanifold}
\begin{equation}
W^{(\mfh)} =\{x\in W\mid \mfs_x\in (\mfh)\}
\end{equation}
which is the union of orbits of the same infinitesimal orbit type $(\mfh)$.  One knows that each connected component of  $W^{(\mfh)}$ is a submanifold. Note that these submanifolds are not $S$-invariant if $S=\Oq$. We call the $S$-orbit of such a submanifold an $S$-component of the orbit type manifold $W^{(\mfh)}$. Moreover we call the stratification on $W$ whose strata are the $S$-components of the orbit type manifolds the {\em infinitesimal} $S$-{\em stratification}.

Let $N$ be an $\oF$-saturated submanifold of $M$, i.e., $\F$-saturated and $\mfa$-invariant. Then $N_P:=\pi^{-1}(N)$ is  $S$-invariant and $\oE$-saturated, i.e., $\E$-saturated and $\mfa$-invariant. Then $N_W:=\rho(N_P)$ is $S$-invariant. Conversely, we may start with an $S$-invariant submanifold $V$ of $W$ to obtain an $\oF$-saturated submanifold $V_M$ of $M$. We thus have a one-to-one correspondence between $\oF$-saturated submanifolds of $M$ and $S$-invariant submanifolds of $W$. In particular $M/\oF=W/S$.

Let $x\in M, \hat x\in P$ and $\bar x\in W$ such that $\pi(\hat x)=x, \rho(\hat x)=\bar x$. Furthermore let $\widehat{L}_{\hat x}$ be the leaf through $\hat x$ and $\overline{\whL}_{\hat x}$ its closure. 

\begin{prop}\label{prop:correspondence}
The $\mfa$-stratification of $(M,\F)$ corresponds to the infinitesimal $S$-stratification of $W$. In fact, for $\hat x\in P, x\in M$ and $\bar x\in W$ with $\pi(\hat x)=x$ and $\rho(\hat x)=\bar x$ we have that $\soq_{\bar x}=\soq_{\overline{\whL}_{\hat x}}=\omega(\mfa_x\hat x)$ is abelian of the same dimension as $\mfa_x$.
\end{prop}
\begin{proof}
We want to show that $\mfa_x$ and $\soq_{\bar x}$ are antiisomorphic as Lie algebras, where $x\in M, \bar x\in W$ and $\hat x\in P$ are points with $\pi(\hat x)=x$ and $\rho(\hat x)=\bar x$. Let $F_{\hat x}=\pi^{-1}(x)\cap \overline{\widehat{L}}_{\hat x}$. Since $F_{\hat x}=\hat x\cdot \SOq_{ \overline{\widehat{L}}_{\hat x}}$ we can view $F_{\hat x}$ as the Lie group $\SOq_{ \overline{\widehat{L}}_{\hat x}}$ by designating $\hat x$ as the neutral element. By Theorem \ref{thm:equivariantMolino} the vector fields on $F_{\hat x}$ induced by $\mfa_x$ are right-$\SOq_{ \overline{\widehat{L}}_{\hat x}}$-invariant and therefore are right-invariant vector fields on the Lie group $\SOq_{ \overline{\widehat{L}}_{\hat x}}$. Since the action of $\mfa_x$ on $F_{\hat x}$ is free and transitive the Lie algebras $\mfa_x$ and $\soq_{\overline{\widehat{L}}_{\hat x}}$ are antiisomorphic. In particular $\soq_{\overline{\widehat{L}}_{\hat x}}$ is abelian.

Let $(X,Y)\in (\mfa\times \soq)_{\hat x}$. Then $0=(X,Y)^*_{\hat x}=X^*_{\hat x}+Y^*_{\hat x}$. Thus $X_{\hat x}^*\in T_{\hat x}\pi^{-1}(x)$ and $Y^*_{\hat x}\in T_{\hat x}\overline{\widehat{L}}_{\hat x}$, so $X\in\mfa_x$ and $Y\in \soq_{\bar x}$.
Since the $\mfa$-action and the $\soq$-action are free, $X$ is determined by $Y$ and vice versa. But for every $X\in\mfa_x$ we have $(X,-\omega(X^*_{\hat x}))\in (\mfa\times \soq)_{\hat x}$, because by definition of a connection form, ${\omega(X^*_{\hat x})}^*_{\hat x}=X^*_{\hat x}$.

Thus the correspondence of the isotropy of the transverse $\mfa$-action on $(M,\F)$ at $x$ and the $\soq$-action on $W$ at $\bar x$ is given by $\omega(\mfa_x\hat x)=\soq_{\overline{\whL}_{\hat x}}=\soq_{\bar x}$.
\end{proof}
In particular, the leaf closure $\oL_x$ of a leaf $L_x$ with isotropy $\mfa_x=\mfa_L=\mfa_{\oL}$ maps to an orbit with isotropy type $(\omega(\mfa_x\hat x))$, where $\hat x\in P$ is an arbitrary point over $x$. 

The proposition implies that the identity components of all isotropy groups $\SOq_{\bar x}$ are abelian. The following lemma gives a little more information in case $\SOq_{\bar x}$ is disconnected.

\begin{lemma} \label{lem:isotropytrivial} For every $\bar x\in W$, the adjoint representation of the isotropy group $\SOq_{\bar x}$ on $\soq_{\bar x}$ is trivial.
\end{lemma}
\begin{proof} By Proposition \ref{prop:correspondence}, every element in $\soq_{\bar x}$ is of the form $\omega(X^*_{\hat x})$ for some $\hat x\in P$ with $\rho(\hat x)=\bar x$ and $X\in \mfa_x$, where $x=\pi(\hat x)$. For $g\in \SOq_{{\widehat{L}}_{\hat x}}$, we have
\begin{equation*}
\begin{aligned}
\Ad_{g} \omega(X^*_{\hat x})&=\omega(X^*_{\hat x} \cdot g^{-1})\\
&=\omega(X^*_{\hat x \cdot g^{-1}}) \\
&=\omega(X^*_{\hat x})
\end{aligned}
\qquad
\begin{aligned}
&\omega\ \mbox{is a connection form}\\
&\mbox{the}\ \mfa\mbox{- and the}\ \SOq\mbox{-action commute} \\
&\hat x \cdot g^{-1}\in \whL_{\hat x}, X^*\ \text{is}\ \E\mbox{-transverse and}\ \omega\ \mbox{is}\ \E \mbox{-basic}\
\end{aligned}
\end{equation*}
By continuity this equation follows for all $g\in \overline{\SOq_{{\widehat{L}}_{\hat x}}}$.  We clearly have $\overline{\SOq_{{\widehat{L}}_{\hat x}}}\subset \SOq_{\overline{\widehat{L}}_{\hat x}}$. For $g\in \SOq_{\overline{\widehat{L}}_{\hat x}}$, there is a sequence $g_n\in \SOq$ satisfying $\hat x g_n\in \whL_{\hat x}$ and $\hat x g_n\to \hat x g$. Because $\SOq$ respects $\E$, we have $g_n\in \SOq_{{\widehat{L}}_{\hat x}}$. As the $\SOq$-action is free, $g_n\to g$. Thus, $\SOq_{\bar x}=\SOq_{\overline{\widehat{L}}_{\hat x}}=\overline{\SOq_{{\widehat{L}}_{\hat x}}}$.
\end{proof}

\begin{prop}\label{prop:lift}
Let $N$ be an $\oF$-saturated submanifold of $M$. Then  we have $\pi^*\nu N\cong \rho^*\nu N_W$.
\end{prop}
\begin{proof}
Let $\cH_1=\cH\subset\nu\E$ be the transverse horizontal distribution of the transversely Riemannian submersion $\pi:P\to M$. Similarly there is a transverse horizontal distribution $\cH_2\subset\nu\E$ for $\rho:P\to W$, defined as the common kernel of the $\theta_i\in \Omega^1(P,\E)$ from Lemma \ref{lem:freetransverseaction}. Then
$$
(TN/T\F)^{\cH_1}_{\hat x}\oplus \hat x\cdot  \soq= T_{\hat x}N_P/T_{\hat x} \E= (TN_W)^{\cH_2}_{\hat x}\oplus \mfa\cdot\hat x,
$$
where the direct sums are orthogonal with respect to the transverse metric $g_P$ and the superscripts denote the respective horizontal lifts, e.g.~$(TN/T\F)^{\cH_1}$ is the geometric realization of $\pi^*(TN/T\F)$ as a subbundle of $\cH_1$. As $\pi$ and $\rho$ are transversely Riemannian submersions, $\pi_*:\cH_1\to \nu \F$ and $\rho_*:\cH_2\to TW$ are pointwise orthogonal, and we have
$$
(\pi^*\nu N)_{\hat x}\cong(\nu N)^{\cH_1}_{\hat x}= \nu_{\hat x} N_P= (\nu N_W)^{\cH_2}_{\hat x}\cong (\rho^*\nu N_W)_{\hat x}.
$$
Here $\nu N=(TN/T\F)^{\perp_g}, \nu N_P= (TN_P/T\E)^{\perp_{g_P}}$ and $\nu N_W=(TN_W)^{\perp_{g_W}}$.
\end{proof}

\begin{rem}\label{rem:isotopystrat}
In \cite{HT2} a similar stratification of $(W,S)$ was considered which we will explain in this remark. We will only use it in Section \ref{sec:polytope}. For a subgroup $H$ of $S$ let $(H)=\{K\mid sHs^{-1}=K\ \text{for some}\ s\in S\}$ be the $S$-conjugacy class of $H$ in $S$. For $x$ in an orbit the conjugacy class $(S_x)$ does not depend on $x$  and is therefore an invariant of $Sx$. The class $(S_x)$ is called the {\em orbit type} of $Sx$. For an isotropy group $H$, define the $(H)$-{\em orbit type submanifold}
$W^{(H)} =\{x\in W\mid S_x\in (H)\}$
which is the union of orbits of the same orbit type $(H)$. We call the stratification of $(W,S)$ whose strata are the $S$-components of these orbit type manifolds the $S$-{\em stratification} (in contrast to the previously defined infinitesimal $S$-stratification). Via the above correspondence it was shown in \cite{HT2} to correspond to the isotopy stratification of $(M,\F)$; the stratum containing a leaf closure $\oL$ is defined to be the union of leaf closures $\overline{L'}$ such that there is an $\F$-isotopy from $L$ to $L'$. The strata of the infinitesimal $S$-stratification of $(M,\F)$ considered before are slightly larger.
\end{rem}

\subsection{Commuting actions}
Let $\F$ be a transversely oriented Killing foliation on a complete manifold $M$ of codimension $q$ with structural Killing algebra $\mfa$. Let $(P,\E)$ be the Molino bundle of oriented orthonormal frames of $\nu \F$. Then we have a natural isometric transverse action of $\mfa\times \soq$ on $(P,\E)$ and we may consider the equivariant cohomology $H^*_{\mfa\times \soq}(P,\E)$. We will show that
\begin{equation}\label{eq:commutingaction}
H^*_\mfa(M,\F)\cong H^*_{\mfa\times\soq}(P,\E)\cong H^*_{\soq}(W),
\end{equation}
where the right hand side is to be understood as the equivariant cohomology of the $\soq^\star$-algebra $\Omega(W)$, see Example \ref{ex:gstaraction}. By Example \ref{ex:weilmodel} it may also be interpreted as the $\SOq$-equivariant cohomology $H^*_\SOq(W)$.
The isomorphisms \eqref{eq:commutingaction} are a special case of the following
\begin{prop}\label{prop:commutingaction}
For every $\oF$-saturated submanifold $N$ in $M$ we have the following isomorphisms of graded rings:
$$
 H^*_\mfa(N,\F)\cong H^*_{\mfa\times\soq}(N_P,\E)\cong H^*_{\soq}(N_W).
$$
\end{prop}
\begin{proof}
Because $\Omega(N_P,\E)_{\bas\, \mfa}=\Omega(N_P,\oE)=\Omega(N_W)$ and $\Omega(N_P,\E)_{\bas\, \soq}=\Omega(N,\F)$, this follows from Proposition \ref{prop:abstractcommutingaction} applied twice to the $(\mfa\times \soq)^\star$-algebra $\Omega(N_P,\E)$ if we can check that $\Omega(N_P,\E)$ is of type (C) both as an $\soq^\star$-algebra and as an $\mfa^\star$-algebra (but of course not as an $(\mfa\times\soq)^\star$-algebra). We have already shown in Lemma \ref{lem:freetransverseaction} that $\Omega(N_P,\E)$  is of type (C) as an $\mfa^\star$-algebra. Let $\alpha_i$ be a basis of $\soq^*$. Then the span of the $\alpha_i\circ\omega$, where $\omega$ is the connection form introduced in Section \ref{sec:Molinobundle} is an $\SOq$-invariant subspace of $\Omega^1(N_P,\E)$ because of the $\SOq$-invariance of $\ker\omega$, i.e., $\Omega(P,\E)$ is of type (C) also as an $\soq^\star$-algebra.
\end{proof}

Although the isomorphism in Proposition \ref{prop:commutingaction} is an isomorphism of graded rings, it also allows to relate the $S(\mfa^*)$-module structure of $H^*_\mfa(M,\F)$ to the $S(\soq^*)^\soq$-module structure of $H^*_\soq(W)$, as we will see in the Proposition \ref{prop:eqformalCM} which will be needed in Sections \ref{sec:basicMorseBott} and \ref{sec:polytope}. See Appendix \ref{app:CMmodules} for the definition and basic properties of graded Cohen-Macaulay rings and modules. For more details, see e.g.~the book \cite{BrunsHerzog}. For general properties of actions whose equivariant cohomology is a Cohen-Macaulay module see \cite{GT2009}. 
\begin{lemma} \label{lem:dimsoqequivcohom}
$\dim_{S(\soq^*)^\SOq} H^*_\SOq(W)$ equals the maximal dimension of an $\SOq$-isotropy algebra.
\end{lemma}
\begin{proof} For torus actions this was proven in \cite[Proposition 5.1]{FrPu}.
Their proof goes through in our situation essentially because all isotropy algebras of the $\SOq$-action are abelian. 

By Lemma \ref{lem:isotropytrivial}, we have $H^*_\SOq(\SOq x)=S(\soq_x^*)^{\SOq_x}=S(\soq_x^*)$ for all $x$. In particular, these rings are Cohen-Macaulay of Krull dimension $\dim \soq_x$ and we have by Proposition \ref{app:propserre} that they are also Cohen-Macaulay as modules over $S(\soq^*)^\soq$ of the same Krull dimension:
$$
\dim_{H^*(B\SOq)}H^*(B\SOq_x)=\dim_{S(\soq^*)^\soq} S(\soq_x^*)=\dim \soq_x.
$$
With this statement proven (instead of \cite[Proposition 3.2]{FrPu}) the rest of  the proof works as in \cite{FrPu}.
\end{proof}
\begin{prop} \label{prop:eqformalCM}
If $M/\oF$ is compact and $\F$ has closed leaves, then the $\mfa$-action on $(M,\F)$ is equivariantly formal if and only if $H^*_\soq(W)$ is a Cohen-Macaulay module over $S(\soq^*)^\soq$.
\end{prop}
\begin{proof}
Proposition \ref{prop:commutingaction} implies that $H^*_\mfa(M,\F)$ is a graded Cohen-Macaulay ring if and only if $H^*_\soq(W)$ is a Cohen-Macaulay ring. Because $H^*_\mfa(M,\F)$ and $H^*_\soq(W)$ are finitely generated modules over $S(\mfa^*)$ respectively $S(\soq^*)^\soq$ (see Corollary \ref{cor:leafspacecptthenfinite} for $H^*_\mfa(M,\F)$) we may apply Proposition \ref{app:propserre} to obtain
\begin{align*}
\dim_{S(\mfa^*)}H^*_\mfa(M,\F)&=\dim_{H^*_\mfa(M,\F)}H^*_\mfa(M,\F)\\
&=\dim_{H^*_\soq(W)}H^*_\soq(W)=\dim_{S(\soq^*)^\soq} H^*_\soq(W)
\end{align*}
and both $H^*_\mfa(M,\F)$ and $H^*_\soq(W)$ are graded Cohen-Macaulay rings if and only if they are graded Cohen-Macaulay modules.

The right hand side of this equation equals the maximal dimension of an $\soq$-isotropy algebra by Lemma \ref{lem:dimsoqequivcohom}. If there is a closed leaf, then this maximal dimension is $\dim \mfa$ by Proposition \ref{prop:correspondence}. But a Cohen-Macaulay module over $S(\mfa^*)$ of Krull dimension $\dim \mfa$ is free, as can be seen from the graded version of the Auslander-Buchs\-baum Theorem \cite[Exercise 19.8]{Eisenbud} and the fact that for a graded module over a polynomial ring, the projective dimension is equal to the length of the minimal free resolution \cite[Corollary 1.8]{Syzygies}. Thus, in case there are closed leaves, $H^*_\mfa(M,\F)$ is a free module over $S(\mfa^*)$ if and only if $H^*_\soq(W)$ is a Cohen-Macaulay module over $S(\soq^*)^\soq$.
By Proposition \ref{prop:freemodequivformal}, freeness of $H^*_\mfa(M,\F)$ is equivalent to equivariant formality of the $\mfa$-action.
\end{proof}

\section{Borel-type localization}\label{sec:Borel}
Let $\F$ be a transversely oriented Killing foliation on a complete manifold $M$ with structural Killing algebra $\mfa$ and let $C$ be the union of closed leaves.  Let $Q(\mfa^*)$ be the field of fractions of $S(\mfa^*)$. For an $S(\mfa^*)$-module $H$ we define the localization of $H$ with respect to the zero ideal of $S(\mfa^*)$ as ${\widehat H}=H\otimes_{S(\mfa^*)}Q(\mfa^*)$. We say $H$ is torsion if for every $h\in H$ there is an element $f\in S(\mfa^*)\backslash\{0\}$ with $fh=0$, which is equivalent to $\widehat H=0$.

Let $N$ be an $\F$-saturated, $\mfa$-invariant submanifold. Then it has an $\F$-saturated, $\mfa$-invariant tubular neighborhood $V$ as can be proven for instance via the correspondence principle: $V:=V'_M$, where $V'$ is an $\SOq$-tubular neighborhood of $N_W$. Note that the inclusion $i:N\to V$ induces an isomorphism on cohomology because the inclusion $i_W:N_W\to V'$ does: the correspondence principle yields the commutative diagram
$$
\xymatrix@R=20pt@C=-4pt{
  H^*_\soq(N_W)  &\cong H^*_{\mfa\times \soq}(N_P,\E)&\cong H^*_\mfa(N,\F)  \\
  H^*_\soq(V')\ar[u]^{{i^*_W}}&\,\cong H^*_{\mfa\times \soq}(V'_P,\E)\ar[u]^{i_P^*}&\,\cong H^*_\mfa(V,\F) \ar[u]^{i^*}.}
$$

\begin{lemma}\label{lem:torsion}
Let $\F$ be a  transversely oriented Killing foliation on a complete manifold $M$ with structural Killing algebra $\mfa$, and let $U\subset M$ be an open, $\F$-saturated subset. If there are no closed leaves in $U$ and if the $\mfa$-action on $(U,\F)$ has only finitely many different isotropy algebras, then $H^*_\mfa(U,\F)$ is torsion.
\end{lemma}
\begin{proof}
We prove the lemma by induction over the maximal length of a chain of inclusions of isotropy algebras. If there is only one isotropy algebra on $(U,\F)$, then the action on $(U,\F)$ is free, and $H^*_\mfa(U,\F)=H^*(U/\oF)$ by Corollary \ref{cor:afree}, which is torsion. Assume we have shown the statement for all open, $\F$-saturated subsets of $M$ with maximal length of a chain of isotropy algebras less than or equal to $n$ and let $U\subset M$ be an open and $\F$-saturated subset whose maximal length is equal to $n+1$. Let $N$ be the set of points with locally maximal isotropy algebra in $U$ and $V$ be a tubular neighborhood of $N$ in $U$. Then in the Mayer-Vietoris sequence (Proposition \ref{prop:mayervietoris})
$$
\ldots\to H^*_\mfa(U,\F)\to H^*_\mfa(V,\F)\oplus H^*_\mfa(U\setminus N,\F)\to H^*_\mfa(V\cap (U\setminus N),\F)\to \ldots
$$
$H^*_\mfa(U\setminus N,\F)$ and $H^*_\mfa(V\cap (U\setminus N),\F)$ are torsion by induction. Let $\mfa_1,\ldots, \mfa_k$ be the isotropy algebras occurring in $N$ and note that $N$ is the disjoint union of the $N^{\mfa_i}=\{p\in N\mid X^*_p=0 \text{ for all } X\in \mfa_i\}$. As $\mfa_i\neq\mfa$ for all $i$,
$$
H^*_\mfa(U,\F)\cong H^*_\mfa(N,\F)=\bigoplus_{i=1}^kH^*_{\mfa}(N^{\mfa_i},\F)=\bigoplus_{i=1}^kH^*(N^{\mfa_i}/\oF)\otimes S(\mfa_i^*)
$$
is torsion, so $H^*_\mfa(M,\F)$ is torsion as well. Here the first equality follows from the observation before the Lemma and the last equality is due to Corollary \ref{cor:oneisotropy} as $N^{\mfa_i}$ is open in the complete manifold $M^{\mfa_i}$.
\end{proof}

\begin{thm} \label{thm:borelloc} Let $\F$ be a  transversely oriented Killing foliation on a complete manifold $M$ with structural Killing algebra $\mfa$ and let $C=M^\mfa$ be the union of closed leaves.
If the transverse action of $\mfa$ has only finitely many different isotropy algebras (e.g.~if $M/\oF$ is compact) then
$$
\widehat H^*_\mfa(M,\F)\to \widehat H^*_\mfa(C,\F)
$$
is an isomorphism.
\end{thm}
\begin{proof}
Let $U$ be an $\F$-saturated, $\mfa$-invariant tubular neighborhood of $C$. Then in the Mayer-Vietoris sequence
$$
\ldots\to H^*_\mfa(M,\F)\to H^*_\mfa(U,\F)\oplus H^*_\mfa(M\setminus C,\F)\to H^*_\mfa(U\cap (M\setminus C),\F)\to \ldots
$$
$H^*_\mfa(M\setminus C,\F)$ and $H^*_\mfa(U\cap (M\setminus C),\F)$ are torsion because of Lemma \ref{lem:torsion}. So $\widehat H^*_\mfa(M,\F)\to \widehat H^*_\mfa(U,\F)\cong \widehat H^*_\mfa(C,\F)$ is an isomorphism.
\end{proof}
\begin{cor}\label{cor:torsionmodule}
Let $(M,\F)$ be as in Theorem \ref{thm:borelloc}. Then the kernel of $H^*_\mfa(M,\F)\to H^*_\mfa(C,\F)$ is the torsion submodule of $H^*_\mfa(M,\F)$.
\end{cor}
\begin{proof}
In the diagram
$$
\xymatrix{
  H^*_\mfa(M,\F) \ar[r]\ar[d] & H^*_\mfa(C,\F)\ar[d]  \\
  {\widehat H}^*_\mfa(M,\F)\ar[r]^\cong & {\widehat H}^*_\mfa(C,\F)}
$$
the bottom horizontal map is an isomorphism by Theorem \ref{thm:borelloc}. As the right vertical map is injective, the kernel of the upper horizontal map equals the kernel of the localization map on the left, which is the torsion submodule.
\end{proof}
The following corollary gives an algebraic criterion for the existence of closed leaves in analogy to the criterion for the existence of fixed points for torus actions given in \cite[Corollary 1, p.~45]{Hsiang}.
\begin{cor}
Let $(M,\F)$ be as in Theorem \ref{thm:borelloc}. Then the following conditions are equivalent:
\begin{enumerate}
\item  $\F$ has a closed leaf.
\item The natural map $S(\mfa^*)\to H^*_\mfa(M,\F)$ defining the $S(\mfa^*)$-module structure on $H^*_\mfa(M,\F)$ is injective.
\item $H^*_\mfa(M,\F)$ is not a torsion module.
\end{enumerate}
\end{cor}
\begin{proof}
Conditions $(2)$ and $(3)$ are clearly equivalent.
If $\F$ does not have a closed leaf, then $C=\emptyset$, so $H^*_\mfa(M,\F)$ is torsion by Corollary \ref{cor:torsionmodule}. If $\F$ has a closed leaf, then $1\in H^*_\mfa(C,\F)=H^*(C,\F)\otimes S(\mfa^*)$ is not a torsion element, so $1\in H^*_\mfa(M,\F)$ is not a torsion element either.
\end{proof}
For the following Theorem, compare \cite[Corollary 2, p.~46]{Hsiang}.
\begin{thm}\label{thm:boreleq} Let $(M,\F)$ be as in Theorem \ref{thm:borelloc}. Then
$$
\dim H^*(C/\F)=\dim H^*(C,\F)\leq \dim H^*(M,\F),
$$
and equality is valid if and only if the $\mfa$-action on $(M,\F)$ is equivariantly formal.
\end{thm}
\begin{proof}
Using the spectral sequence in Theorem \ref{thm:specseq} we have
\begin{align*}
\dim H^*(M,\F)&=\dim_{Q(\mfa^*)}\widehat E_1
\geq \dim_{Q(\mfa^*)}\widehat E_\infty \\
&=\dim_{Q(\mfa^*)} \widehat H^*_\mfa(M,\F)=\dim_{Q(\mfa^*)} \widehat H^*_\mfa(C,\F)=\dim H^*(C/\F),
\end{align*}
where the second to last equation is due to Theorem \ref{thm:borelloc}. Equality holds if and only if the spectral sequence collapses at the $E_1$-term, which by Proposition \ref{prop:E2collapse} holds if and only if the action is equivariantly formal.
\end{proof}

\begin{prop}\label{prop:bottomstrclosed}
Let $(M,\F)$ be as in Theorem \ref{thm:borelloc}. If the $\mfa$-action is equivariantly formal, then the bottom stratum of the $\mfa$-action is the union $C$ of the closed leaves.
\end{prop}
\begin{proof}
Consider the diagram
$$
\xymatrix{
  H^*_\mfa(M,\F) \ar[r]\ar[d]^\cong & H^*_\mfa(C,\F)\ar[d]^\cong  \\
  H^*_\soq(W)\ar[r]& H^*_\soq(C_W)}
$$
where the horizontal maps are the natural restriction maps and the vertical isomorphisms are given by the correspondence principle. The upper horizontal map is injective by Corollary \ref{cor:torsionmodule}, so it follows that $H^*_\soq(W)\to H^*_\soq(C_W)$ is injective. If there existed a component $N$ of the bottom stratum of the $\mfa$-action not consisting of closed leaves, then $N_W$ would be a component of the bottom stratum of the $\soq$-action with empty intersection with $C_W$. Let $(\mfh)$ be the infinitesimal isotropy type of $N_W$. By Proposition \ref{prop:correspondence}, $\mfh$ is abelian. Because $N_W$ is a component of the bottom stratum, $\mfh$ acts without fixed vectors on $\left. \nu N_W\right|_{N_W^\mfh}$. Then the composition of the equivariant Thom homomorphism of $N_W\hookrightarrow W$ with the restriction to $N_W$,
$$
H^*_\soq(N_W)\to H^*_\soq(W)\to H^*_\soq(N_W),
$$ is multiplication with the equivariant Euler class of the normal bundle $\nu N_W$, which by Proposition \ref{prop:eulernozero} is not a zero-divisor. Note that $\nu N_W$ is orientable by Proposition \ref{prop:lift} as $\nu N$ carries an orientation induced by the $\mfa$-action. Thus, the Thom homomorphism $H^*_\soq(N_W)\to H^*_\soq(W)$ is injective. This contradicts the fact that the image of $1\in H^0_\soq(N_W)$ under the Thom homomorphism is mapped to zero when restricted to $C_W$.
\end{proof}

\section{(Equivariant) basic Morse-Bott theory}\label{sec:basicMorseBott}
We start by recalling the main notions of Morse-Bott theory. Let $f:M\to\RR$ be a Morse-Bott function on a manifold $M$ and let $N$ be a non-degenerate critical manifold of $f$ in $M$. It follows by the Morse-Bott Lemma that the critical submanifolds are isolated. We have an orthogonal decomposition of the normal bundle $\nu N=\nu^+N\oplus\nu^-N$ into subbundles according to the positive/negative eigenvalues of the Hessian of $f$ on $\nu N$. The fibre dimension $\lambda_N$ of $\nu^-N$ is the {\em index} of $N$ as a critical manifold of $f$. Assuming that the negative normal bundle of each $N$ is orientable, we count a non-degenerate manifold $N$ of $f$ with the polynomial
$$
M_t(f,N):=\sum t^i\dim H^i_c(\nu^-N)\stackrel{{\rm Thom}}{=}t^{\lambda_N}P_t(N),
$$
where $P_t(N)=\sum_it^i\dim H^i(N)$ is the Poincar\'e polynomial of $N$.
Then we define
$$
M_t(f):=\sum_N M_t(f,N).
$$
In this generality, the appearing coefficents might be infinite. 
Assume $M$ is a complete manifold and $\F$ is a transversely oriented Killing foliation on $M$ with Killing algebra $\mfa$ such that $M/\oF$ is compact. 
The basic Poincar\'e polynomial $P_t(M,\F)$ and the equivariant basic Poincar\'e series $P_t^\mfa(M,\F)$ are defined by
$$
P_t(M,\F)=\sum_it^i\dim H^i(M,\F)\qquad \mbox{and}\qquad P_t^\mfa(M,\F)=\sum_it^i\dim H^i_\mfa(M,\F).
$$
Let $f:M\to\mR$ be an $\F$-basic Morse-Bott function. Since it is constant on the leaves of $\F$ it is also constant on the leaves of $\oF$ and therefore $\mfa$-invariant. The (equivariant) basic Morse polynomial/series are defined by
$$
M_t(f,\F):=\sum_N t^{\lambda_N}P_t(N,\F)\qquad \mbox{and}\qquad M_t^\mfa(f,\F):=\sum_N t^{\lambda_N}P_t^\mfa(N,\F),
$$
where $N$ runs over the critical manifolds of $f$. All appearing coefficients in these polynomials and series are finite by Proposition \ref{prop:basiccohomfinite} and Corollary \ref{cor:leafspacecptthenfinite}.

The function $f_P=f\circ \pi$ is $\mfa\times \soq$-invariant and $\E$-basic. This function descends to an $\SOq$-invariant function $f_W:W\to\mR$. If $N$ is a critical submanifold of $f$ then $N_W$ is a critical submanifold of $f_W$, and conversely every critical submanifold of $f_W$ arises this way. We have $\pi^*\Hess f=\Hess f_P=\rho^*\Hess f_W$, and hence under the isomorphism of Proposition \ref{prop:lift} the subbundle $\pi^*\nu^- N\subset \pi^*\nu N$ maps to the subbundle $\rho^*\nu^- N_W\subset \rho^*\nu N_W$, and the indices of $N$ and $N_W$ coincide, i.e., $\lambda_N=\lambda_{N_W}$. Furthermore, by Proposition \ref{prop:lift} the negative bundle of $N$ is orientable if and only if that of $N_W$ is. In this case, the $\SOq$-equivariant Morse polynomial of $f_W$ is
$$
M_t^\SOq(f_W)=\sum_Nt^{\lambda_N}P_t^{\SOq}(N_W),
$$
where $N$ runs over the critical manifolds of $f$. Recall that the equivariant Poincar\'e series is defined as
$$
P_t^\SOq(N)=\sum_it^i\dim H_\SOq^i(N).
$$
Besides the well-known notion of an (equivariantly) perfect Morse-Bott function, we also have the natural notion of an (equivariantly) perfect basic Morse-Bott function.
\begin{prop}\label{prop:Morse}
Let $M$ be a complete manifold with a transversely oriented Killing foliation $\F$ such that $M/\oF$ is compact. Then
$$
P_t^\mfa(M,\F)=P_t^{\SOq}(W).
$$
Let $f:M\to\mR$ be an $\F$-basic Morse-Bott function such that the negative normal bundles of the critical submanifolds of $f$ are orientable. Then
$$
M_t^\mfa(f,\F)=M_t^{\SOq}(f_W).
$$
In particular, $f$ is an equivariantly perfect basic Morse-Bott function if and only if $f_W$ is equivariantly perfect.
\end{prop}
\begin{proof}
The first equality is due to Proposition \ref{prop:commutingaction}, because $P_t^\mfa(M,\F)$ and $P_t^{\SOq}(W)$ are the Poincar\'e series of $H^*_\mfa(M,\F)$ and $H^*_\soq(W)$, respectively.
For any connected component $N$ of the critical set of $f$, we have $H^*_\mfa(N,\F)\cong H^*_\soq(N_W)$ by Proposition \ref{prop:commutingaction} and $\lambda_N=\lambda_{N_W}$. Thus
$$
M_t^\mfa(f,\F)\stackrel{\mbox{\tiny def}}{=}\sum_N t^{\lambda_N}P_t^\mfa(N,\F)=\sum_N t^{\lambda_{N_W}}P_t^\soq(N_W)\stackrel{\mbox{\tiny def}}{=}M_t^{\SOq}(f_W).
$$
\end{proof}
As a corollary we obtain equivariant basic Morse inequalities.
\begin{thm}\label{thm:basicMorseBott}
Let $M$ be a complete manifold with a transversely oriented Killing foliation $\F$ such that $M/\oF$ is compact. For an $\F$-basic Morse-Bott function $f:M\to\mR$ such that the negative normal bundles of the critical submanifolds of $f$ are orientable we have the equivariant basic Morse-Bott inequalities
$$
M_t^\mfa(f,\F)-P_t^\mfa(M,\F)\geq (1+t)Q(t)
$$
for some polynomial $Q(t)$ with nonnegative coefficients.
\end{thm}
Recall that $P_t\geq Q_t$ for polynomials/series means that $P_t-Q_t$ has only nonnegative coefficients.
\begin{proof}
We have $M_t^\soq(f)-P_t^\soq(W)\geq (1+t)Q(t)$ by the classical equivariant Morse-Bott theory \cite{AtiyahBott1982}. Therefore the result follows from Proposition \ref{prop:Morse}.\end{proof}

\begin{thm}\label{thm:morsebott}
Let $\F$ be a transversely oriented Killing foliation on a complete manifold $M$ with structural Killing algebra $\mfa$ such that $M/\oF$ is compact.  Let $f:M\to\RR$ be a basic Morse-Bott function such that its critical set is the union $C$ of the closed leaves. Then $f$ is $\mfa$-equivariantly perfect and the $\mfa$-action on $(M,\F)$ is equivariantly formal.
\end{thm}
\begin{proof}
The $\mfa$-action on the negative normal bundle has no fixed vectors and hence yields an orientation on the negative normal bundles of the components of $C$. Thus, the negative normal bundles of the $N_W$ are orientable. Since the critical set of $f_W$ is the bottom stratum (see Section \ref{sec:stratification}) of the $\SOq$-action, this function is $\SOq$-equivariantly perfect by Theorem \ref{thm:equivperfect}. By Proposition \ref{prop:Morse}, $f$ is equivariantly perfect, i.e., $P_t^\mfa(M,\F)=M_t^\mfa(f,\F)$.

One shows as in \cite[Theorem 7.1]{GT2009} that the existence of the Morse-Bott function $f_W$ with critical set equal to the set of points with isotropy algebra of dimension $\dim \mfa$ implies that $H^*_\SOq(W)$ is Cohen-Macaulay of Krull dimension $\dim\mfa$: For a real number $a$, let $W^a=f^{-1}((-\infty,a])$. For a critical value $\kappa$, denote by $N^\kappa_i$ the connected components of the critical set at level $\kappa$ and write $\lambda^\kappa_i:=\lambda_{N^\kappa_i}$. Using Proposition \ref{prop:eulernozero} as in  \cite{GT2009}, the long exact cohomology sequence of the pair $(W^{\kappa+\varepsilon},W^{\kappa-\varepsilon})$ splits into a short exact sequence of $S(\soq^*)^\soq$-modules
\[
0\to H^*_\soq(W^{\kappa+\varepsilon},W^{\kappa-\varepsilon}) \to H^*_\soq(W^{\kappa+\varepsilon})\to H^*_\soq(W^{\kappa-\varepsilon})\to 0.
\]
To apply induction we need to show that if the module on the right is Cohen-Macaulay of Krull dimension $\dim \mfa$, then so is the one in the middle. By Lemma \ref{lem:sesCM} it suffices to show that the module on the left is Cohen-Macaulay of the same Krull dimension. By the Thom isomorphism $\bigoplus_i H^*_\soq((N^\kappa_i)_W)\to H^*_\soq(W^{\kappa+\varepsilon},W^{\kappa-\varepsilon})$ it suffices to show that $\bigoplus_i H^*_\soq((N^\kappa_i)_W)$ is Cohen-Macaulay as a module (or equivalently as a ring, by Proposition \ref{app:propserre}). But by the correspondence principle, see Proposition \ref{prop:commutingaction}, this ring is isomorphic to the Cohen-Macaulay ring $\bigoplus_i H^*_\mfa(N^\kappa_i)=\bigoplus_i H^*(N^\kappa_i/\oF)\otimes S(\mfa^*)$.

Having proven that $H^*_\SOq(W)$ is Cohen-Macaulay, Proposition  \ref{prop:eqformalCM} implies that the $\mfa$-action is equivariantly formal.
\end{proof}

\begin{thm}\label{thm:basicperfect}
Under the assumptions of Theorem \ref{thm:morsebott}, $f$ is a perfect basic Morse-Bott function, i.e.
$$
P_t(M,\F)=\sum_N t^{\lambda_N} P_t(N/\F),
$$
where $N$ runs over the components of $C$.
\end{thm}
\begin{proof}
Since $P_t(S(\mfa^*))P_t(M,\F)= P_t^\mfa(M,\F)$ and $P_t^\mfa(M,\F)=M_t^\mfa(f,\F)$ by Theorem \ref{thm:morsebott}, we may calculate
\begin{align*}
P_t(S(\mfa^*))P_t(M,\F)&= P_t^\mfa(M,\F)\\
&=M_t^\mfa(f,\F)=\sum_Nt^{\lambda_N}P_t^\mfa(N,\F)=P_t(S(\mfa^*))\sum_Nt^{\lambda_N}P_t(N,\F).
\end{align*}
Cancelling $P_t(S(\mfa^*))$, the desired equality follows.
\end{proof}
\begin{rem}
We proved Morse inequalities for equivariant basic cohomology in Theorem \ref{thm:basicMorseBott}. For (nonequivariant) basic cohomology, they were proven in \cite{Alvarez}. Theorem \ref{thm:basicperfect} gives a sufficient condition for the latter to be equalities.

Another lower bound for the total number of critical leaf closures of a basic function (not necessarily Morse-Bott) was given in \cite{HT2} in terms of the essential LS-category.
\end{rem}
\begin{ex}\label{ex:hopf}
In Example \ref{ex:hopfexample}, it is easy to find an $\F$-basic Morse-Bott function $f:S^3\to\RR$ whose critical set is the union of the two compact leaves. This is a perfect basic Morse-Bott function with one critical leaf of index $0$ and one of index $2$. This shows that the basic Betti numbers $b_i(M,\F)$ are $1$ for $i=0,2$ and $0$ otherwise.
\end{ex}

\section{The basic Poincar\'e series}\label{sec:Poincareseries}
\subsection{The $G$-equivariant Gysin sequence}\label{sec:Gysin}
In this section $W$ will denote an arbitrary $G$-manifold, where $G$ is a compact connected Lie group of rank $r$. Later we will apply the result in this section to the $\SOq$-action on the space $W=P/\oE$, where $P$ is the Molino bundle associated to a transversely oriented Killing foliation $(M,\F)$. For $i=1,\ldots, r$ define
$$
W_i:=\{x\in W\mid \rk \mfg_{x}\geq i\}
$$
and
$$
W_{(i)}:=\{x\in W\mid \rk \mfg_{x}= i\}.
$$
For every point $x\in W$, the isotropy algebras of nearby points are conjugate to subalgebras of $\mfg_x$, so the function $x\mapsto \rk \mfg_x$ is upper semi-continuous. Thus, $W_{(i)}$ is closed in $W\backslash W_{i+1}$. We will assume that the action satisfies the following condition:
\begin{equation}\label{condition}
\mfg_x\subsetneq\mfg_y\Longrightarrow \rk \mfg_x<\rk\mfg_y
\end{equation}
for all $x,y\in W$. This is fulfilled for instance if all isotropy algebras are abelian, which is the case for the $\SOq$-action on $W=P/\oE$ by Proposition \ref{prop:correspondence}.
The condition implies that the components of the $W_{(i)}$ coincide with the components of the infinitesimal orbit type manifolds $W^{(\mfh)}=\{x\in W\mid \mfg_x\in (\mfh)\}$; here $(\mfh)$ is the set of conjugacy classes of $\mfh$ in $\mfg$.

In analogy to Theorem 2 of \cite{Duflot} we obtain the following theorem whose proof may be found in the appendix (Theorem \ref{thmapp:poincareseries}).
\begin{thm}\label{thm:poincareseries}
Consider an action of a compact connected Lie group $G$ on a manifold $W$ satisfying Condition \eqref{condition}, such that all normal bundles $\nu Y$ of the components $Y$ of the infinitesimal orbit type manifolds are orientable. Then
$$
P_t^G(W)=\sum_{Y} t^{\codim Y} P_t^G(Y),
$$
where $Y$ runs over the components of the infinitesimal orbit type manifolds.
\end{thm}
Note that we allow the coefficients of the Poincar\'e series to be infinite. The sums over $Y$ are well-defined as all appearing summands are nonnegative. A sufficient condition for finiteness of the coefficients is $\dim H^*(W)<\infty$.

In the language of Kirwan \cite[p.~33]{Kirwan} the theorem says that the stratification of $W$ by components of $\mfg$-orbit type manifolds is equivariantly perfect.

\subsection{Poincar\'e series of $H^*_\mfa(M,\F)$}
Let $(M,\F)$ be a Killing foliation with structural Killing algebra $\mfa$. Let
$$
M_i:=\{x\in M\mid \rk \mfa_{x}\geq i\}
$$
and
$$
M_{(i)}:=\{x\in M\mid \rk \mfa_{x}= i\}.
$$
The latter is the union of leaves $L$ whose {\em transverse dimension} ${\rm trdim\,}L:=\dim\oL-\dim L=\dim\mfa-\dim\mfa_L$ is equal to $\dim\mfa-i$. Under the correspondence of Section \ref{sec:stratification} the $M_{(i)}$ respectively $M_i$ map to $W_{(i)}$ respectively $W_i$ as defined in Section \ref{sec:Gysin}. More precisely, by Proposition \ref{prop:correspondence} a component $X$ of $M_{(i)}$ with $i$-dimensional isotropy algebra $\mfa_X$ maps to a component $X_W$ of $W_{(i)}$ with isotropy type $(\omega(\mfa_X\cdot\hat x))$, where $\hat x\in P$ is some point with $\pi(\hat x)\in X$.

\begin{thm}\label{thm:basicpoincareseries}
Let $(M,\F)$ be a Killing foliation with structural Killing algebra $\mfa$. Then
$$
P_t^\mfa(M,\F)=\sum_{X} t^{\codim X}P_t^\mfa(X,\F)=\sum_{X} \frac{t^{\codim X}}{(1-t^2)^{\dim \mfa_X}}P_t(X/\oF),
$$
where $X$ runs over all components of infinitesimal orbit type manifolds $M^{(\mathfrak b)}$.
\end{thm}
\begin{rem}
The term $\frac{1}{(1-t^2)}$ is to be understood as the series $\sum_{i=0}^\infty t^{2i}$. Thus the equation makes sense even in case a coefficient in some  $P_t(X/\oF)$ is infinity.
\end{rem}
\begin{proof} The $\mathfrak b$-action on the normal bundle $\nu X$ of a component $X$ of $M^{(\mathfrak b)}$ induces an orientation on $\nu X$. By Proposition \ref{prop:lift}, it follows that $\nu X_W$ is orientable. Thus we may calculate, using Theorem \ref{thm:poincareseries} applied to the $\SOq$-action on $W=P/\oE$:
\begin{align*}
P_t^\mfa(M,\F) = P_t^\SOq(W)=\sum_{X} t^{\codim X_W} P_t^\SOq(X_W)=\sum_{X} t^{\codim X}P_t^\mfa(X,\F).
\end{align*}
Here we applied Proposition \ref{prop:commutingaction} on $M$ and all $X$, and we used the fact that $\codim X=\codim X_W$ for all $X$ by Proposition \ref{prop:lift}.
For the second equality in the statement of the theorem, note that $H^*_\mfa(X,\F)=H^*(X/\oF)\otimes S(\mfa_X^*)$ by Corollary \ref{cor:oneisotropy}, so
$$
P_t^\mfa(X,\F) = P_t(S(\mfa_X^*))P_t(X/\oF)=\frac{1}{(1-t^2)^{\dim\mfa_X}}P_t(X/\oF).
$$
\end{proof}

\begin{cor}\label{cor:basicps} Assume that the $\mfa$-action on $(M,\F)$ has only finitely many isotropy algebras, that each orbit type manifold has only finitely many connected components, and that each such component $X$ satisfies $\dim H^*(X/\oF)<\infty$. (These conditions are satisfied if $M/\oF$ is compact.)
If the $\mfa$-action on $(M,\F)$ is equivariantly formal, then
$$
P_t(M,\F)=\sum_{X} {t^{\codim X }}{(1-t^2)^{{\rm trdim\,}L_X}}P_t(X/\oF),
$$
where $L_X$ is any leaf in $X$.
\end{cor}
\begin{proof}
By definition of equivariant formality, $P_t^\mfa(M,\F)=P_t(S(\mfa^*))P_t(M,\F)=\frac{P_t(M,\F)}{(1-t^2)^{\dim\mfa}}$. Since all involved coefficients are finite by assumption, we may multiply the equality in Theorem \ref{thm:basicpoincareseries} by $(1-t^2)^{\dim\mfa}$.  The result follows because of ${\rm trdim\,}L_X=\dim\mfa-\dim\mfa_X$.
\end{proof}

This means that in order to compute the basic Betti numbers in the equivariantly formal case it suffices to understand the stratification of $M/\oF$ with the respective isotropy dimensions.

\begin{ex} In the situation of Example \ref{ex:torusaction}, assume that there is a closed $B$-orbit. As the $T$-orbits are given by the closures of the $B$-orbits, the minimal dimension of a $T$-orbit is $\dim B$. Thus, the maximal isotropy dimension of the $T$-action is $\dim \mft/\mfb=\dim \mfa$. By \cite[Proposition 5.1]{FrPu} the $S(\mft^*)$-module $H^*_T(M)$ is of Krull dimension $\dim \mfa$. The calculation in Example \ref{ex:torusaction} shows that $H^*_T(M)$ is a Cohen-Macaulay module over $S(\mft^*)$ if and only if $H^*_\mfa(M,\F)$ is Cohen-Macaulay over $S(\mfa^*)$. Thus, by Proposition \ref{prop:freemodequivformal}, the $T$-action on $M$ is Cohen-Macaulay in the sense of \cite{GT2009} if and only if the $\mfa$-action on $(M,\F)$ is equivariantly formal.
\end{ex}

Under the finiteness assumptions of Corollary \ref{cor:basicps} the Euler characteristic $\chi(M,\F)=P_{-1}(M,\F)$ of $H^*(M,\F)$ is well-defined, and the corollary implies
\begin{cor} \label{cor:eulercharacteristic}
If the $\mfa$-action is equivariantly formal, then
$$
\chi(M,\F)=\chi(C/\F),
$$
where $C$ is the union of closed leaves.
\end{cor}

\begin{proof}
The summands of the right hand side of the equality in Corollary \ref{cor:basicps} with positive transverse dimension vanish after inserting $-1$. Furthermore the codimension of a component of $C$ is always even, since there is an $\mfa$-representation without trivial summand on its normal space.
\end{proof}

\begin{rem} By the classical Poincar\'{e}-Hopf Theorem, we have the equality of classical Euler characteristics $\chi(M)=\chi(C)$ for any Killing foliation without the assumption of equivariant formality, because a generic element in $\mfa$ yields a vector field on $M$ whose zero set is $C$.

There is a basic version of the Hopf Index Theorem \cite[Theorem 3.18]{BeEfRi} which requires the closed leaves to be isolated. Still without requiring equivariant formality of the $\mfa$-action, this theorem, applied to a generic element in $\mfa$, implies the equality in Corollary \ref{cor:eulercharacteristic}.
\end{rem}

\section{Killing foliations whose space of leaf closures is a polytope}\label{sec:polytope}
In this section we will give another sufficient condition for equivariant formality. Let $\F$ be a transversely oriented Killing foliation of codimension $q$ on a complete manifold $M$ with structural Killing algebra $\mfa$.

\begin{prop}\label{prop:qdima}
If there exists a closed leaf then $2\dim \mfa\leq q$. If there is an isolated closed leaf then $q$ is even.
\end{prop}
\begin{proof}
The isotropy type of a closed leaf is $\mfa$. There is an isotropy algebra of the $\SOq$-action on $W$ antiisomorphic to $\mfa$ by Proposition \ref{prop:correspondence}. But a maximal abelian subalgebra of $\soq$ is $\left[q/2\right]$-dimensional.
Now let $L$ be an isolated closed leaf. Then the $\mfa$-action on the $q$-dimensional normal space $\nu_pL$ decomposes $\nu_p L$ into its two-dimensional weight spaces.
\end{proof}

We will now assume that $M/\oF$ is homeomorphic to a simple, convex polytope in such a way that the projection of the isotopy stratification of $(M,\F)$ introduced in Remark \ref{rem:isotopystrat} to $M/\oF$ coincides with the natural stratification by the faces. For short, we say that $M/\oF$ is a simple, convex polytope. Note that it follows that the $\mfa$-action has only finitely many isotropy algebras.

\begin{thm}\label{thm:polytope}
Assume $M/\oF$ is a simple, convex polytope of dimension $\dim\mfa$ such that the vertices correspond to the closed leaves. Then the $\mfa$-action is equivariantly formal.
\end{thm}
\begin{proof}
We will show that the $S$-stratification on $W$ coincides with the infinitesimal $S$-stratification.
Because of
$$
q- \dim\mfa=\dim M-\dim\F-\dim\mfa=\dim M/\oF=\dim\mfa
$$
we have $q=2\dim\mfa$. The $\mfa$-representation on the normal space $\nu_x L$ of a closed leaf $L$ has $\dim\mfa$ linearly independent weights, as it otherwise would not be effective. This shows that a neighborhood of a closed leaf in $M/\oF$ looks like $(\RR^+)^{\dim\mfa}$ with respect to the infinitesimal $S$-stratification. In a simple polytope there are exactly $n \choose i$ faces of dimension $i$ adjacent to a vertex. As each stratum of the (non-infinitesimal) $S$-stratification is contained in a stratum of the infinitesimal $S$-stratification, it follows that the two stratifications are identical.

The proof for equivariant formality is adapted from \cite[Theorem 3.1]{davisjan}. Let $\tau:M\to M/\oF$ be the projection. A choice of a generic linear functional $\phi$ on the polytope $M/\oF$ defines a decomposition of $M/\oF$ as follows: Let $L_i$ be the closed leaves ordered via $\phi\circ\tau$ with $L_0$ corresponding to the absolute minimum. Define $M_0:=L_0$. If $M_i$ (not to be confused with the $M_i$ from Section \ref{sec:Poincareseries}) is already defined, let $F_{i+1}$ be the closed face defined by all edges at $\tau(L_{i+1})$ below  $\tau(L_{i+1})$  with respect to $\phi$, or equivalently all edges connecting $\tau(L_{i+1})$ with $\tau(M_i)$. Let $M_{i+1}=\tau^{-1}(F_{i+1})\cup M_{i}$.
Denote the corresponding spaces  in $W$ by $W_{i}=(M_{i})_W$ (see Section \ref{sec:stratification}) and $(L_i)_W=\SOq p_i$ for some $p_i\in W$.  We will show by induction that for all $i$, $H^*_\SOq(W_{i})$ is a Cohen-Macaulay module over $S(\soq^*)^\soq$ of Krull dimension $\dim \mfa$. For $i=0$ this is clear as $$
H^*_\SOq(W_{0})=H^*(B\SOq _{p_0})=S(\soq^*_{p_0})^{\SOq_{p_0}^*}=S(\soq^*_{p_0}).
$$ Here we used the standard result about the equivariant cohomology of a homogeneous space \cite[p.~34]{Hsiang} and Lemma \ref{lem:isotropytrivial}. See also the proof of Lemma \ref{lem:dimsoqequivcohom}.

Assume that $H^*_\SOq(W_i)$ is Cohen-Macaulay. We have that $W_{i+1}\setminus \SOq p_{i+1}$ is $\SOq$-homotopy equivalent to $W_i$, as can be seen e.g.~by successively moving down the skeleta: start by lifting a deformation of the (contractible) interior of the top-dimensional face in $W_{i+1}\setminus W_i$ into a tubular neighborhood of the lower skeleton which is $\SOq$-homotopy equivalent to the lower skeleton itself, etc. Now Proposition \ref{prop:Gysinsubm} yields the following short exact sequence:
$$
0\to H^*_\SOq(\SOq p_{i+1})\to H^*_\SOq (W_{i+1})\to H^*_\SOq(W_{i})\to 0.
$$
As $H^*_\SOq(\SOq p_{i+1})$ and $H^*_\SOq(W_{i})$ are Cohen-Macaulay modules of the same Krull dimension $\dim \mfa$ (the same argument as above for $H^*_\SOq(W_0)$ applies), $H^*_\SOq (W_{i+1})$ is Cohen-Macaulay of Krull dimension $\dim \mfa$ as well, see Lemma \ref{lem:sesCM}. Thus, $H^*_\soq(W)$ is a Cohen-Macaulay module, and by Proposition \ref{prop:eqformalCM} it follows that the $\mfa$-action is equivariantly formal.
\end{proof}

As a consequence of Theorem \ref{thm:polytope} we have a formula relating the basic Betti numbers and the number of faces of $M/\oF$. First, we need a lemma.

\begin{lemma} \label{lem:codimX}
If $M/\oF$ is a simple, convex polytope of dimension $\dim \mfa$ such that the vertices correspond to the closed leaves, then for each component $X$ of an orbit type manifold, we have $\codim X=2\dim \mfa_X$. Moreover, for every leaf $L\subset X$, we have $\trdim L=\dim X/\oF$.
\end{lemma}
\begin{proof}
By Theorem \ref{thm:polytope}, the formula holds for $X$ any closed leaf. For arbitrary $X$, choose a closed leaf $L\subset \overline{X}$, and some point $x\in L$. The isotropy representation of $\mfa$ on the $2\dim\mfa$-dimensional normal space $\nu_xL$ has exactly $\dim\mfa$ weights $\alpha_1,\ldots,\alpha_{\dim \mfa}$, which hence are linearly independent. Write $\nu_xL=\bigoplus_{i=1}^{\dim \mfa} V_{\alpha_i}$, where $V_{\alpha_i}$ is the two-dimensional weight space corresponding to the weight $\alpha_i$. As $X$ is a component of the orbit type manifold corresponding to $\mfa_X$, we have $T_x\overline{X}=T_xL\oplus \bigoplus_{i:\, \alpha_i|_{\mfa_X}=0}V_{\alpha_i}$ and hence $\nu_x\overline{X}=\bigoplus_{i:\, \alpha_i|_{\mfa_X}\neq 0}V_{\alpha_i}$. As the weights are linearly independent, there are exactly $\dim \mfa_X$ weights that do not vanish on $\mfa_X$. We thus have $\codim X=2\dim \mfa_X$ and  $\dim X/\oF=\dim \mfa-\dim \mfa_X=\trdim L$.
\end{proof}

\begin{thm} \label{thm:polytopebasicps} If $M/\oF$ is a simple, convex polytope of dimension $\dim \mfa$ such that the vertices correspond to the closed leaves. Then 
$$
P_t(M,\F)=\sum_i \lambda_i t^{q-2i}(1-t^2)^i,
$$
where $\lambda_i$ denotes the number of faces of $M/\oF$ of dimension $i$.
\end{thm}
\begin{proof} The $\mfa$-action is equivariantly formal by Theorem \ref{thm:polytope}. Hence, by Corollary \ref{cor:basicps},
$$
P_t(M,\F)=\sum_X t^{\codim X}(1-t^2)^{\trdim L_X},
$$
where $X$ runs over the components of infinitesimal orbit type manifolds. If $X/\oF$ is an open face of the polytope $M/\oF$ of dimension $i$, then by Lemma \ref{lem:codimX}, we have $\trdim L_X=\dim X/\oF=i$ and $\codim X=2\dim \mfa_X=2\dim \mfa-2\trdim L_X=q-2i$.\end{proof}

\begin{ex} Let us recalculate Example \ref{ex:hopfexample} with the new methods. $S^3/\oF=S^2/S^1$ is a line segment, with the end points corresponding to the compact leaves. Thus, in the notation above, $\lambda_0=2$ and $\lambda_1=1$. Because $q=2$, it follows $P_t(S^3,\F)=1+t^2$.
\end{ex}

\appendix
\section{Some results on $G$-equivariant cohomology} \label{sec:appA}
\begin{lemma}\label{lem:quotient}
Let $G$ be a compact connected Lie group acting on a manifold $M$ and $H$ be a connected normal closed Lie subgroup acting trivially. Then
$$
H_G^*(M)=H_{G/H}^*(M)\otimes S(\mfh^*)^\mfh.
$$
\end{lemma}
\begin{proof}
Let $\mfk$ be an ideal in the Lie algebra $\mfg$ of $G$ complementary to $\mfh$, i.e., $\mfg=\mfk\oplus\mfh$ as Lie algebras. Then
$$
(S(\mfg^*)\otimes\Omega(M))^\mfg=(S(\mfh^*)\otimes S(\mfk^*)\otimes \Omega(M))^{\mfh\oplus\mfk}=S(\mfh^*)^\mfh\otimes (S(\mfk^*)\otimes \Omega(M))^\mfk,
$$
because $\mfh$ acts trivially on $S(\mfk^*)\otimes\Omega(M)$. Since $d_\mfg(f\otimes \omega)=f\otimes d_\mfk(\omega)$ for any $f\in S(\mfh^*)^\mfh$ and $\omega\in S(\mfk^*)\otimes\Omega(M)$ we obtain
$$
H^*_G(M)=H^*_\mfg(M)=S(\mfh^*)^\mfh \otimes H_\mfk^*(M)=S(\mfh^*)^\mfh \otimes H_{G/H}^*(M)
$$
as $S(\mfh^*)^\mfh\otimes S(\mfk^*)^\mfk=S(\mfg^*)^\mfg$-algebras using \eqref{eq:inf}.
\end{proof}

\begin{lemma}\label{lem:normalizer}
Let $G$ be a compact Lie group acting on a manifold $M$ with exactly one infinitesimal isotropy type $(\mfh)$. Then
$$
H^*_G(M)=H^*_{N_G(\mfh)}(M^\mfh)
$$
where $N_G(\mfh)$ is the normalizer of $\mfh$ in $G$.
\end{lemma}
\begin{proof}
Since $G$ acts with one infinitesimal isotropy type, we have $M=G\times_{N_G(\mfh)}M^\mfh$. Then
$$
M_G=EG\times_G M=(EG\times_G G)\times_{N_G(\mfh)}M^\mfh=EG\times_{N_G(\mfh)}M^\mfh=(M^\mfh)_{N_G(\mfh)},
$$
which implies the lemma.
\end{proof}

Recall that for a $G$-oriented $G$-equivariant vector bundle $E\to M$ of rank $m$  the equivariant Euler class $e_G(E)\in H^m_G(M)$ is defined as the classical Euler class of the vector bundle $E_G\to M_G$. See e.g.~\cite[Section C.6.2]{GGK}.
\begin{lemma}\label{lem:reduction}
Let $G$ be a compact Lie group acting on a manifold $M$ with one infinitesimal isotropy type $(\mfh)$. Furthermore let $E\to M$ be a $G$-oriented $G$-equivariant vector bundle. Under the isomorphism $H^*_G(M)=H^*_{N_G(\mfh)}(M^\mfh)$ of Lemma \ref{lem:normalizer} we have the equality of equivariant Euler classes
$$
e_G(E)=e_{N_G(\mfh)}(E|_{M^\mfh}).
$$
\end{lemma}
\begin{proof}
The horizontal maps in the following diagram define an isomorphism of vector bundles:
$$
\xymatrix{
 E_G=EG \times_G E \ar[r] \ar[d] & EG\times_{N_G(\mfh)} E|_{M^\mfh} =(E|_{M^\mfh})_{N_G(\mfh)} \ar[d]  \\
 M_G=EG\times_G M \ar[r] & EG\times_{N_G(\mfh)} M^\mfh  = (M^\mfh)_{N_G(\mfh)}
}
$$
Therefore the Euler classes coincide.
\end{proof}
\begin{lemma}\label{lem:restricteuler}
Let $G$ be a compact Lie group and $H\subset G$ a closed subgroup. Furthermore let $E\to M$ be a  $G$-oriented $G$-vector bundle. Then under the natural map $H^*_G(M)\to H^*_{H}(M)$, the equivariant Euler class $e_G(E)$ maps to $e_H(E)$.
\end{lemma}
\begin{proof} Let $f:M_H=EG \times_H M\to EG\times_G M=M_G$ be the natural projection. The statement follows from the isomorphism of $H$-equivariant vector bundles $E_H=f^*(E_G)$.
\end{proof}
The next proposition generalizes \cite[Proposition 4]{Duflot} from actions of tori to actions of compact connected Lie groups.
\begin{prop}\label{prop:reduction}
Let $G$ be a compact Lie group acting on a manifold $M$, such that the connected component $H$ of the kernel of the action is positive dimensional. Furthermore let $E\to M$ be a $G$-oriented $G$-vector bundle of rank $2m$ on which a maximal torus $T$ of $H$ acts without fixed vectors. Then the $G$-equivariant Euler class $e_G(E)$ is not a zero-divisor in $H^*_G(M)$.
\end{prop}
\begin{proof}
Denote by $G^0$ the connected component of $G$. As $H^*_G(M)=H^*(EG \times_G M)=H^*_{G^0}(M)^{G/G^0}$ by \cite[Cor.~III.2.3]{Borel}, the natural map $H^*_G(M)\to H^*_{G^0}(M)$ is injective. Thus by Lemma \ref{lem:restricteuler} if we can show that $e_{G^0}(E)\in H^*_{G^0}(E)$ is not a zero-divisor, then $e_G(E)\in H^*_G(M)$ is not a zero-divisor as well.

As $G^0$ leaves invariant the connected components of $M$, and $e_{G^0}(M)$ consists of the equivariant Euler classes of the bundle $E$ restricted to the components, we can without loss of generality assume that $M$ is connected. Under the isomorphism $H^*_{G^0}(M)=H^*_{G^0/H}(M)\otimes S(\mfh^*)^\mfh$ of Lemma \ref{lem:quotient} we consider the equivariant Euler class $e_{G^0}(E)\in H^{2m}_{G^0}(M)$ as a polynomial in $\mfh$ with values in $H^*_{G^0/H}(M)$. To show that $e_{G^0}(E)$ is not a zero-divisor it is sufficient to show that the part with polynomial degree $m$, i.e., its part in $H^0_{G^0/H}(M)\otimes S^m(\mfh^*)^\mfh$ does not vanish. Again by Lemma \ref{lem:restricteuler} in the following diagram $e_{G^0}$ restricts to $e_H$, which in turn restricts to $e_T$.
$$
\xymatrix@R=20pt@C=-4pt{
 e_{G^0}(E)\in \ar@{|->}[d]&H^*_{G^0}(M)  \ar[d] &=H^*_{G^0/H}(M)\otimes S(\mfh^*)^\mfh  \ar[d]  \\
 e_H(E)\in \ar@{|->}[d]&H^*_H(M)  \ar[d]&=\,\, H^*(M)\otimes S(\mfh^*)^\mfh  \ar[d]\\
 e_T(E)\in &H^*_T(M)  & =\,\,\, H^*(M)\otimes S(\mft^*)
}
$$
On the right hand side $H^0_{G^0/H}(M)\otimes S^m(\mfh^*)^\mfh$ maps injectively to $H^0(M)\otimes S^m(\mft^*)$.  Hence the relevant portion of $e_{G^0}(E)$ is determined by the respective part of $e_T(E)$. As $T$ acts without fixed vectors on $E$ this is, up to a scalar, the product of the $T$-weights on $E$ (\cite[Examples C.41/C.44]{GGK}), in particular nonzero.
\end{proof}

\begin{prop}\label{prop:eulernozero}
Let $G$ be a compact Lie group acting on a manifold $M$ with exactly one infinitesimal isotropy type $(\mfh)$ and let $E\to M$ be a $G$-oriented $G$-vector bundle. Furthermore assume that a maximal abelian subalgebra of $\mfh$ acts without fixed vectors on $E|_{M^\mfh}$. Then the equivariant Euler class $e_G(E)$ is not a zero-divisor of $H^*_G(M)$.
\end{prop}
\begin{proof}
By Lemma \ref{lem:reduction} it is enough to consider the $N_G(\mfh)$-action on $M^\mfh$. Then the statement follows from Proposition \ref{prop:reduction}.
\end{proof}
\begin{rem}
The condition that a maximal abelian subalgebra of $\mfh$ acts without fixed vectors cannot be replaced by the weaker condition of $\mfh$ itself acting without fixed vectors. Consider any representation of a compact Lie group $G$ without fixed vectors such that a maximal torus $T$ has fixed vectors, e.g.~the adjoint representation of any semisimple Lie group. Then $e_T=0\in H^*_T({\rm pt})$ and since $H^*_G({\rm pt})\to H^*_T({\rm pt})$ is injective, $e_G=0$.
\end{rem}

From now on, we restrict to $G$-actions satisfying the condition
\begin{equation}\label{app:condition}
\mfg_x\subsetneq\mfg_y\Longrightarrow \rk \mfg_x<\rk\mfg_y.
\end{equation}
As in Section \ref{sec:Gysin}, we will for $i=1,\ldots, r$ denote
$$
W_i:=\{x\in W\mid \rk \mfg_{x}\geq i\}
$$
and
$$
W_{(i)}:=\{x\in W\mid \rk \mfg_{x}= i\}.
$$
Recall from Section \ref{sec:Gysin} that this condition implies that the components of $W_{(i)}$ coincide with the components of the infinitesimal orbit type manifolds $W^{(\mfh)}$. The next theorem is well-known for (Hamiltonian) torus actions, see, e.g., \cite[Section I.5]{Kirwan} or \cite[Proposition 2.1]{TolmanWeitsman}.

\begin{thm}\label{thm:equivperfect}
Consider an action of a compact Lie group $G$ on a manifold $M$ satisfying Condition \eqref{app:condition}. If $f:M\to\RR$ is a $G$-invariant Morse-Bott function such that its critical set is equal to a union of connected components of the bottom stratum, and such that the negative normal bundles of the critical manifolds are $G$-orientable, then $f$ is equivariantly perfect.
\end{thm}
\begin{proof}
Atiyah-Bott show in \cite[Proposition 1.9]{AtiyahBott1982} that the relative cocycles are self-completing if we can show that the respective equivariant Euler-classes of the negative normal bundles are not zero-divisors.
 If for some critical manifold $Y$ with isotropy type $(\mfh)$ there was a normal vector in the restricted negative normal bundle $\nu^- Y|_{Y^{\mfh}}$ fixed by some maximal abelian subalgebra of $\mfh$, then there would exist a point close to $Y$ whose isotropy algebra has rank equal to $\rk\mfh$. This would contradict  Condition \eqref{app:condition}.
Thus we may apply Proposition \ref{prop:eulernozero} to the negative normal bundles of the components of the critical manifolds.
\end{proof}

\begin{thm}\label{thmapp:poincareseries}
Consider an action of a compact Lie group $G$ on a manifold $W$ satisfying Condition \eqref{app:condition} such that all normal bundles $\nu Y$ of the $G$-components of the infinitesimal orbit type manifolds are $G$-orientable. Then
$$
P_t^G(W)=\sum_{Y} t^{\codim Y} P_t^G(Y).
$$
\end{thm}

\begin{proof}
Consider for any component $Y$ of $W_{(i)}$ the normal bundle $\nu Y$ of $Y$ in $W\setminus W_{i+1}$. As in the proof of Theorem \ref{thm:equivperfect}, Condition \eqref{app:condition} implies the assumptions of Proposition \ref{prop:eulernozero} with $M=Y$ and $E=\nu Y$ are satisfied, so $e_G(\nu Y)$ is not a zero-divisor of $H^*_G(Y)$. Thus the equivariant push-forward $H^*_G(Y)\to H^*_G(W\setminus W_{i+1})$ of the inclusion $Y\subset W\setminus W_{i+1}$ is injective because composition with the restriction to $Y$ is multiplication with $e_G(Y)$. Note that the push-forward increases the degree by $\codim Y$. The same argument as in the proof of Theorem 1 in \cite{Duflot} shows that $e_G(W_{(i)})$ is not a zero-divisor in $H^*_G(W_{(i)})$ and the equivariant push-forward   $H^*_G(W_{(i)})\to H^*_G(W\setminus W_{i+1})$ of the inclusion $W_{(i)}\hookrightarrow W\setminus W_{i+1}$ is injective. On a summand $H^*_G(Y)$ of $H^*_G(W_{(i)})$ the push-forward shifts the degree by $\codim Y$.
We have shown that the Gysin sequence (\cite[p.~255]{Duflot}) splits into short exact sequences
\begin{equation}\label{eqn:Gysin}
0\to H^*_G(W_{(i)})\to H^*_G(W-W_{i+1}) \to H^*_G(W-W_i)\to 0.
\end{equation}
Now the statement follows from a dimension count using \eqref{eqn:Gysin} inductively just as in the proof of Theorem 2 in \cite{Duflot}. 
\end{proof}
On the way we have shown:
\begin{prop}\label{prop:Gysinsubm}
Consider an action of a compact Lie group $G$ on a manifold $W$ satisfying Condition \eqref{app:condition}. If $N$ is a closed invariant submanifold on which $G$ acts with only one infinitesimal isotropy type such that the normal bundle $\nu N$ of $N$ is $G$-orientable, then the Gysin sequence splits into short exact sequences
$$
0\to H^*_G(N)\to H^*_G(W) \to H^*_G(W-N)\to 0.
$$
\end{prop}

\section{Cohen-Macaulay modules} \label{app:CMmodules}

Using the language of e.g.~\cite[Section 1.5]{BrunsHerzog}, a graded ring $R$ (graded over the integers) is *local if it has a unique *maximal ideal, where a *maximal ideal is a graded ideal $\mfm\neq R$ which is maximal among the graded ideals. Thus, for a compact connected Lie group $G$ with Lie algebra $\mfg$, $S(\mfg^*)^G$ is a Noetherian graded *local ring. Note that in general a *maximal ideal is not necessarily maximal, but for $S(\mfg^*)^\mfg$ this is the case.

Let $R$ be a Noetherian graded *local ring, with *maximal ideal $\mfm$. Then the {\it depth} of a finitely generated graded module $A$ over $R$ is defined as the length of a maximal $A$-regular sequence in $\mfm$:
\[
\depth A=\grade(\mfm,A).
\]
The {\it Krull dimension} of $A$, denoted $\dim_R A$ or simply $\dim A$ is defined as the Krull dimension of the ring $R/\Ann(A)$, where $\Ann(A)=\{r\in R\mid rA=0\}$, i.e., the supremum of the lengths of chains of prime ideals in $R$ containing $\Ann(A)$. We always have $\depth A\leq \dim A$ and we define

\begin{defn} A finitely generated graded module $A$ over a Noetherian graded *local ring $R$ is  {\emph Cohen-Macaulay} if $\depth A=\dim A$. The ring $R$ is said to be {\emph Cohen-Macaulay} if $R$ is Cohen-Macaulay as a module over itself.
\end{defn}
The following lemma is standard. See e.g.~\cite[Lemma 5.4]{GT2009}.
\begin{lemma}\label{lem:sesCM}
Let $0\to A\to B\to C\to 0$ be an exact sequence of finitely generated graded modules over a Noetherian graded *local ring $R$.  If $A$ and $C$ are Cohen-Macaulay of the same Krull dimension $n$, then $B$ is also Cohen-Macaulay of Krull dimension $n$.
\end{lemma}

Let $\mfg$ be a finite-dimensional Lie algebra and $A$ a $\mfg^\star$-algebra. Then the equivariant cohomology $H^*_\mfg(A)$ of $A$  as defined in \eqref{eq:cartandef} is a graded algebra over $S(\mfg^*)^\mfg$. Graded versions of \cite[Propositions IV.B.11, IV.B.12]{Serre} imply \begin{prop}\label{app:propserre}
If $H^*_\mfg(A)$ is finitely generated as a module over $S(\mfg^*)^\mfg$, then
$
\dim_{S(\mfg^*)^\mfg} H^*_\mfg(A)=\dim_{H^*_\mfg(A)}H^*_\mfg(A)$ and $\depth_{S(\mfg^*)^\mfg} H^*_\mfg(A)=\depth_{H^*_\mfg(A)}H^*_\mfg(A)$. In particular, $H^*_\mfg(A)$ is a Cohen-Macaulay module over $S(\mfg^*)^\mfg$ if and only if it is a Cohen-Macaulay ring. 
\end{prop}
We will apply this proposition to two different types of $\mfg^\star$-algebras: the equivariant cohomology $H^*_G(M)$ of an action of a compact connected Lie group $G$ on a manifold $M$, and the equivariant basic cohomology $H^*_\mfa(M,\F)$ of the transverse action of the structural Killing algebra $\mfa$ of a Killing foliation $\F$ on a complete manifold $M$. In the first case, the assumption on finite generation is satisfied if $\dim H^*(M)<\infty$ by \cite[Theorem 6.6.1]{GS1999}, and in the second case if $M/\oF$ is compact, see Proposition \ref{cor:leafspacecptthenfinite}.



\begin{thebibliography}{10}

\bibitem[AP 1993]{allday}
C.~Allday and V.~Puppe, {\em Cohomological methods in transformation groups}, Cambridge Studies in Advanced Mathematics, 32. Cambridge University Press, Cambridge, 1993.

\bibitem[Alv 1993]{Alvarez}
J.~Alvarez L\'opez,
{\em Morse inequalities for pseudogroups of local isometries},
J.~Differential Geom.~{\bf 37} (1993), no.~3, 603--638.

\bibitem[AB 1982]{AtiyahBott1982}
M.~Atiyah and R.~Bott, {\em The Yang-Mills equations over Riemann surfaces}, Philos. Trans.~Roy.~Soc.~London Ser.~A  {\bf 308}  (1983),  no.~1505, 523--615.

\bibitem[BPR 2003]{BeEfRi}
V.~Belfi, E.~Park and K.~Richardson,
{\em A Hopf index theorem for foliations},
Differential Geom.~Appl.~{\bf 18} (2003), no. 3, 319--341.

\bibitem[BBFMP 1960]{Borel} A.~Borel, G.~Bredon, E.~Floyd, P.~Montgomery, and R.~Palais, {\em Seminar on transformation groups}, Annals of Math.~Studies 46,  Princeton University Press, Princeton, 1960.

\bibitem[Bot 1967]{Bott} R.~Bott, {\em Vector fields and characteristic numbers}, Michigan Math.~J.~{\bf 14} (1967), 231--244.



\bibitem[BH 1993]{BrunsHerzog} W.~Bruns and J.~Herzog, {\em Cohen-Macaulay rings}, Cambridge Studies in Advanced Mathematics, 39. Cambridge University Press, Cambridge, 1993.

\bibitem[Car 1949]{cartan} H.~Cartan, {\em Cohomologie r\'eelle d'un espace fibre\'e principal differentiable}, S\'eminaire Cartan, expos\'es 19 et 20, 1949.

\bibitem[Car 1950]{cartan2}
H.~Cartan, {\em Notions d'alg\`{e}bre diff\'{e}rentielle; application aux groupes de Lie et aux vari\'{e}t\'{e}s o\`{u} op\`{e}re un groupe de Lie}, Colloque de topologie, Bruxelles, 1950, 15--27.

\bibitem[DJ 1991]{davisjan}
M.~W.~Davis and T.~Januszkiewicz, {\em Convex polytopes, Coxeter orbifolds and torus actions},  Duke Math.~J.~{\bf 62}  (1991),  no.~2, 417--451.

\bibitem[Duf 1983]{Duflot} J.~Duflot, {\em Smooth toral actions}, Topology {\bf 22} (1983), 253--265.


\bibitem[Eis 1995]{Eisenbud} D.~Eisenbud, {\em Commutative algebra with a view toward algebraic geometry}, Graduate Texts in Mathematics, 150. Springer-Verlag, New York, 1995.

\bibitem[Eis 2004]{Syzygies} D.~Eisenbud, {\em The geometry of syzygies}, Graduate Texts in Mathematics, 229. Springer-Verlag, New York, 2004.


\bibitem[ESH 1985]{AlSeHe} A.~El Kacimi-Alaoui, V.~Sergiescu and G.~Hector, {\em La cohomologie basique d'un feuilletage riemannien est de dimension finie}, Math.~ Z.~{\bf 188}  (1985),  no.~4, 593--599.


\bibitem[FP 2003]{FrPu} M.~Franz and V.~Puppe, {\em Exact sequences for equivariantly formal spaces}, preprint. \texttt{arXiv:math/0307112}.

\bibitem[Ghy 1984]{Ghys}
{{\'E.} Ghys},
{\em Feuilletages riemanniens sur les vari\'et\'es simplement connexes},
 Ann.~Inst.~Fourier (Grenoble) {\bf 34} (1984), no.~4, 203--223.

\bibitem[GT 2009]{GT2009}
{O.~Goertsches and D.~T\"oben},
{\em Torus actions whose equivariant cohomology is Cohen-Macaulay}, preprint. 
{\texttt{arXiv:0912.0637}.}

\bibitem[GKM 1998]{GKM} M.~Goresky, R.~Kottwitz, and R.~MacPherson, {\em Equivariant cohomology, Koszul duality, and the localization theorem}, Invent.~Math.~{\bf 131} (1998), no.~1, 25--83.


\bibitem[GGK 2002]{GGK}
V.~Guillemin, V.~Ginzburg, and Y.~Karshon, {\em Moment maps, cobordisms, and Hamiltonian group actions}, Mathematical Surveys and Monographs, 96. AMS, Providence, 2002.


\bibitem[GLS 1996]{GLS1996}
{V.~Guillemin, E.~Lerman, S.~Sternberg},
{\em Symplectic fibrations and multiplicity diagrams}, Cambridge University Press, Cambridge, 1996.


\bibitem[GS 1999]{GS1999}
{V.~Guillemin, S.~Sternberg},
{\em Supersymmetry and Equivariant de Rham Theory},
Springer-Verlag, Berlin, 1999.


\bibitem[Hsi 1975]{Hsiang}
 W.-Y.~Hsiang, {\em Cohomology theory of topological transformation groups}, Ergebnisse der Mathematik und ihrer Grenzgebiete, Band 85, Springer-Verlag, New York-Heidelberg, 1975.

\bibitem[HT 2009]{HT2}
{S.~Hurder and D.~T\"{o}ben},
{\em Transverse LS category for Riemannian foliations}, Trans.~Amer. Math.~Soc.~{\bf 361} (2009), no.~11, 5647--5680.

\bibitem[KT 1975]{kambertondeur}
F.~Kamber and P.~Tondeur, {\em Foliated bundles and characteristic classes}, Lecture Notes in Mathematics, Vol.~493, Springer-Verlag, Berlin-New York, 1975.

\bibitem[Kir 1984]{Kirwan}
{F.~Kirwan},
{\em Cohomology of quotients in symplectic and algebraic geometry},
Mathematical Notes, 31. Princeton University Press, Princeton, 1984.

\bibitem[McC 2001]{McC}
{J.~McCleary}, {\em A user's guide to spectral sequences}, second edition, Cambridge Studies in Advanced Mathematics, 58. Cambridge University Press, Cambridge, 2001.

\bibitem[Mol 1988]{Molino}
{P.~Molino},
{\em Riemannian foliations},
 with appendices by G.~Cairns, Y.~Carri\`ere, \'E.~Ghys, E.~Salem and V.~Sergiescu,
Birkh\"auser Boston Inc., Boston, 1988.

\bibitem[Moz 1985]{Mozgawa}
W.~Mozgawa, {\em Feuilletages de Killing},
Collect.~Math.~{\bf 36} (1985), no.~3, 285--290.

\bibitem[Rei 1959]{Reinhart} B.~Reinhart, {\em Harmonic integrals on foliated manifolds}, Amer.~J.~Math.~{\bf 81} (1959), 529--536.


\bibitem[Sal 1988]{Salem1988}
{{\'E}.~Salem},
{\em Une g\'en\'eralisation du th\'eor\`eme de Myers-Steenrod aux pseudogroupes d'isom\'e\-tries},
 Ann.~Inst.~Fourier (Grenoble) {\bf 38} (1988), no.~2, 185--200.

\bibitem[Sat 1956]{Satake1956}
I.~Satake, 
{\em On a generalization of the notion of manifold},
Proc.~Nat.~Acad.~Sci.~U.S.A.~{\bf 42} (1956), 359--363.

\bibitem[Ser 2000]{Serre} J.-P.~Serre, {\em Local algebra}. Springer Monographs in Mathematics. Springer-Verlag, Berlin, 2000.

\bibitem[Ste 1974]{Stefan} P.~Stefan, {\em
Accessibility and foliations with singularities},
Bull.~Amer.~Math.~Soc.~{\bf 80} (1974), 1142--1145.

\bibitem[Sus 1973]{Sussmann} H.~Sussmann, {\em Orbits of families of vector fields and integrability of distributions},  Trans. Amer.~Math.~Soc.~{\bf 180}  (1973), 171--188.

\bibitem[TW 1999]{TolmanWeitsman} S.~Tolman and J.~Weitsman, {\em On the cohomology rings of Hamiltonian $T$-spaces}, Northern California Symplectic Geometry Seminar, Amer.~Math.~Soc.~Transl.~Ser.~2, {\bf 196} (1999), 251-258.


\end{thebibliography}
\end{document}